\documentclass[a4paper,12pt]{amsart}
\usepackage{tcolorbox}
\usepackage{pdflscape}
\usepackage[justification=centering,labelformat=empty]{caption}
\usepackage{amsopn,amsfonts,amssymb,amsthm,mathtools,amsmath}
\usepackage{color}
\usepackage{dsfont}
\usepackage{latexsym}

\usepackage{wrapfig}
\usepackage[T1]{fontenc}
\usepackage{lmodern}
\usepackage{textcomp}
\usepackage{soul} 
\usepackage{multicol}

\usepackage{hyperref}
\hypersetup{colorlinks=true,linkcolor=blue,citecolor=red}
\usepackage[utf8]{inputenc}
\usepackage{float}
\usepackage{hhline}
\usepackage{multirow}
\usepackage{anysize}
\marginsize{1.4cm}{1.4cm}{1.4cm}{1.4cm}
\usepackage{graphicx}

\newcommand{\cabe}[1]{\textcolor{green}{#1}}

\def\b{\mathfrak{b}}

\def\g{\mathfrak{g}}
\def\h{\mathfrak{h}}

\def\a{\mathfrak{a}}
\def\s{\mathfrak{s}}

\def\m{\mathfrak{m}}
\def\l{\mathfrak{l}}

\def\hcx{\{J_{\alpha}\}}

\def\su{\mathfrak{su}}

\def\gl{\mathfrak{gl}}
\def\sl{\mathfrak{sl}}
\def\hol{\mathfrak{hol}}

\def\X{\mathfrak{X}}

\def\B{\mathcal{B}}

\def\C{\mathbb{C}}
\def\R{\mathbb{R}}

\def\N{\mathbb{N}}

\def\H{\mathbb{H}}

\newcommand{\GL}{\operatorname{GL}}
\newcommand{\SL}{\operatorname{SL}}

\def\ad{\operatorname{ad}}

\def\Id{\operatorname{Id}}

\def\I{\operatorname{I}}

\newcommand{\End}{\operatorname{End}}

\newcommand{\ncp}{\nabla^\mathrm{CP}} 
\newcommand{\nob}{\nabla^\mathrm{Ob}}

\def\alt{\raise1pt\hbox{$\bigwedge$}}

\def\la{\langle}
\def\ra{\rangle}
\def\multiset#1#2{\ensuremath{\left(\kern-.3em\left(\genfrac{}{}{0pt}{}{#1}{#2}\right)\kern-.3em\right)}}

\theoremstyle{plain}
\newtheorem{theorem}{\bf Theorem}[section]
\newtheorem{corollary}[theorem]{\bf Corollary}
\newtheorem{proposition}[theorem]{\bf Proposition}
\newtheorem{lemma}[theorem]{\bf Lemma}

\theoremstyle{definition}

\theoremstyle{remark}
\newtheorem{remark}[theorem]{\bf Remark}

\newcommand{\ri}{{\rm (i)}}
\newcommand{\rii}{{\rm (ii)}}

\newcommand{\matriz}[1]{\ensuremath{\begin{bmatrix}#1\end{bmatrix}}}

\newcommand{\tr}{\operatorname{tr}}

\newenvironment{smallarray}[1]
{\null\,\vcenter\bgroup\scriptsize
	\renewcommand{\arraystretch}{0.7}%
	\arraycolsep=.13885em
	\hbox\bgroup$\array{@{}#1@{}}}
{\endarray$\egroup\egroup\,\null}

\title{Hypercomplex structures on special linear groups}

\author{Adrián Andrada}
\email{adrian.andrada@unc.edu.ar}

\author{Agustín Garrone}
\email{agustin.garrone@mi.unc.edu.ar}

\author{Alejandro Tolcachier}
\email{atolcachier@unc.edu.ar}

\address{FAMAF, Universidad Nacional de C\'ordoba and CIEM-CONICET, Av. Medina Allende s/n, X5000HUA C\'ordoba, Argentina}

\subjclass[2020]{53C15, 53C26, 53C29.}
\keywords{Hypercomplex manifold, holonomy, Obata connection, Lie group, special linear group.}

\begin{document}
	
	\begin{abstract}
		The purpose of this article is twofold. First, we prove that the 8-dimensional Lie group $\SL(3,\R)$ does not admit a left-invariant hypercomplex structure. To accomplish this we revise the classification of left-invariant complex structures on $\SL(3,\R)$ due to Sasaki. Second, we exhibit a left-invariant hypercomplex structure on $\SL(2n+1,\C)$, which arises from a complex product structure on $\SL(2n+1,\R)$, for all $n\in \N$. We then show that there are no HKT metrics compatible with this hypercomplex structure. Additionally, we determine the associated Obata connection and we compute explicitly its holonomy group, thus providing a new example of an Obata holonomy group properly contained in $\GL(m,\mathbb{H})$ and not contained in $\SL(m,\mathbb{H})$, where $4m=\dim_\R \SL(2n+1,\C)$. 
	\end{abstract}

	\maketitle
	
	\section{Introduction}
	Hypercomplex manifolds are quaternionic analogues of complex manifolds. Indeed, a smooth manifold $M$ is called hypercomplex if it admits a hypercomplex structure, that is, a triple of complex structures $\{J_1,J_2,J_3\}$ obeying the laws of the quaternions (see equation \eqref{eq: quat}). In particular, for all $p\in M$, the tangent space $T_p M$ admits a quaternionic action, thus the dimension of the manifold is $4n$ for some $n\in \N$. 
	
	Boyer has classified in \cite{Boyer} compact hypercomplex manifolds in real dimension 4. They are 4-dimensional tori, quaternionic Hopf surfaces and $\mathrm{K3}$ surfaces. However, the classification of compact hypercomplex manifolds of dimension 8 is far from complete. 
	
	An important family of hypercomplex manifolds is given by hyper-K\"ahler manifolds, that is, hypercomplex manifolds equipped with a Riemannian metric $g$, Hermitian with respect to each complex structure (i.e., a hyper-Hermitian metric), such that $(J_\alpha, g)$ is K\"ahler for all $\alpha=1,2,3$. Hyper-K\"ahler manifolds have Riemannian holonomy in $\operatorname{Sp}(n)$, which implies that they are Ricci-flat. However, the class of hypercomplex manifolds is much broader than the one of hyper-K\"ahler manifolds.
	
	Another valuable source of examples of hypercomplex manifolds is provided by compact Lie groups. Indeed, Joyce proved in \cite{Joyce} the remarkable fact that any compact Lie group becomes hypercomplex (with a left-invariant hypercomplex structure) after it is multiplied by a sufficiently large torus. 
	Hypercomplex structures have also been studied on other families of Lie groups. For instance, there are many examples of nilmanifolds (i.e. compact quotients of nilpotent Lie groups by discrete cocompact subgroups) admitting left-invariant hypercomplex structures (see \cite{BDV,DF2,DF,DF3}). Since non-toral nilmanifolds do not admit K\"ahler structures (see \cite{BG}), these hypercomplex nilmanifolds do not admit any hyper-K\"ahler metric. More recently, almost abelian Lie groups (and their quotients by co-compact discrete subgroups) equipped with left-invariant hypercomplex structures were studied (see \cite{AB, AB1, LT}). 
	
	M.~Obata showed in \cite{Ob} that any $4n$-dimensional hypercomplex manifold admits a unique torsion-free connection $\nob$ such that each complex structure $J_\alpha$, $\alpha=1,2,3$, is parallel. This connection is called the \textit{Obata connection} and its holonomy group, $\operatorname{Hol}^{\mathrm{Ob}}$, is contained in the quaternionic linear group $\GL(n,\H)$. Note that in the case of a hyper-K\"ahler manifold the Obata connection coincides with the Levi-Civita connection of the hyper-K\"ahler metric, and therefore the Obata (or Riemannian) holonomy is contained in $\operatorname{Sp}(n)$. Another distinguished subgroup of $\operatorname{GL}(n, \H)$ is its commutator subgroup $\operatorname{SL}(n, \H)$; we point out that $\operatorname{Sp}(n)$, $\operatorname{SL}(n,\H)$ and $\operatorname{GL}(n,\H)$ itself are the only subgroups of $\operatorname{GL}(n,\H)$ that appear in the Merkulov-Schwachh\"ofer list of possible irreducible holonomy groups of a torsion-free linear connection (see \cite{MeS}). $4n$-dimensional hypercomplex manifolds whose Obata holonomy is contained in $\SL(n,\H)$ are called nowadays $\SL(n,\H)$-manifolds, and such a reduction of the holonomy imposes some restrictions on the underlying manifold. For instance, it was proved by Verbitsky \cite[Claim 1.2]{Ver} that if $(M,\hcx)$ is an $\SL(n,\H)$-manifold then the canonical bundle of the complex manifold $(M,J_\alpha)$ is holomorphically trivial for any $\alpha=1,2,3$. The study of hypercomplex manifolds with Obata holonomy contained in $\operatorname{SL}(n,\H)$ is very active; see for instance \cite{GLV,IP,LW,LW1,SV}.
	
	Even though the Obata holonomy group is an important invariant of hypercomplex manifolds, it is not often determined explicitly. In fact, examples where the Obata holonomy has been explicitly computed are scarce. For instance, it was shown in \cite{Sol} that the Obata holonomy on $\operatorname{SU}(3)$, equipped with the left-invariant hypercomplex structure constructed by Joyce in \cite{Joyce}, coincides with $\operatorname{GL}(2,\H)$. On the other hand, according to \cite{BDV}, the Obata holonomy of a hypercomplex nilmanifold  always lies in $\operatorname{SL}(n,\H)$. In the solvmanifold setting, it was shown in \cite{AB} that the Obata connection of a hypercomplex almost abelian solvmanifold is always flat, even though its holonomy may not be contained in $\operatorname{SL}(n,\H)$.
	
	A class of hyper-Hermitian metrics less restrictive than hyper-K\"ahler metrics are the so-called hyper-K\"ahler with torsion (or HKT) metrics \cite{HP}. These are hyper-Hermitian structures satisfying $\partial\Omega=0$, where $\Omega=\omega_2+i\omega_3$ and $\partial$ is the Dolbeault differential of $(M,J_1)$. It is known (see \cite{GP}) that a hyper-Hermitian structure $(\hcx, g)$ on $M$ is HKT if and only if $\nabla^b_1=\nabla^b_2=\nabla^b_3$, where $\nabla^b_\alpha$ is the Bismut connection associated with $(J_\alpha ,g), \; \alpha =1,2,3$. For instance, all the hypercomplex structures constructed by Joyce on compact Lie groups admit HKT metrics.
	We note that the class of hyper-K\"ahler manifolds is strictly contained in the class of HKT manifolds (see \cite{DF}) and this, in turn, is strictly contained in the class of hypercomplex manifolds (see \cite{FG}). The geometry of HKT manifolds has been intensively studied in recent years (see \cite{BDV, BGV, FuGe, GeTa, GLV, IP, LW}). 
	
	The purpose of this article is to study the existence of left-invariant hypercomplex structures on some special linear groups. First, we address the natural question of whether the 8-dimensional real simple Lie group $\SL(3,\R)$ admits a left-invariant hypercomplex structure. By exploiting the classification of left-invariant complex structures on this Lie group given by T. Sasaki in \cite{Sasaki} we are able to show the non-existence of such a hypercomplex structure (see Theorem \ref{thm: SL(3,R) no tiene HCS}). Furthermore, we exhibit a hypercomplex structure on $\SL(3,\R)$ which is not left-invariant; this construction holds in fact for $\SL(2n+1,\R)$ for arbitrary $n$, and essentially follows from the Iwasawa decomposition of $\SL(2n+1,\R)$ and the results of Joyce \cite{Joyce} (see Remark \ref{obs: there is a hypercomplex structure}). We point out that we needed to carry out a revision of Sasaki's result, exhibiting in particular a new equivalence between two complex structures previously believed to be non-equivalent. 
	
	Secondly, we focus on the complex special linear groups $\SL(2n+1,\C)$. By applying results in \cite{AS}, we construct a left-invariant hypercomplex structure on $\SL(2n+1,\C)$ for all $n\in \N$ (see Corollary \ref{cor: hpcx underlying}). Moreover, we show that this hypercomplex structure does not admit any compatible left-invariant HKT metric (Proposition \ref{prop: there is no compatible HKT metric}). Finally, in Section \ref{sec: holonomy} we study the holonomy group of the Obata connection associated to this hypercomplex structure on $\SL(2n+1,\C)$, using some new results proved in Proposition \ref{prop: complex-hol} and Corollary \ref{coro: holonomies} regarding the Obata holonomy of hypercomplex structures arising from the method given in \cite{AS}. In particular, by computing explicitly the Obata holonomy algebra we show that $\SL(2n+1,\C)$ is not an $\SL(m,\H)$-manifold, where $4m=\dim_\R \SL(2n+1,\C)$. Remarkably, this Obata holonomy algebra appears to be new in the existing literature, and in particular, it is not contained in the Merkulov-Schwachhöfer list. Note that this is not in contradiction with the aforementioned classification, as the holonomy representation in question is reducible.  
	
	
	\smallskip 
	
	\textbf{Acknowledgments.} 
	This work was partially supported by CONICET, SECyT-UNC and FONCyT (Argentina). The authors would like to thank Giovanni Gentili and the anonymous referee for many helpful comments. 
	
	\medskip
	
	\section{Preliminaries} \label{section: preliminares}
	\subsection{Complex structures} 
	An almost complex structure on a differentiable manifold $M$ is an automorphism
	$J$ of the tangent bundle $TM$ satisfying $J^2=-\Id $, where $\Id$ is the identity endomorphism of $TM$. Certainly, the existence of an almost complex structure on $M$ forces the dimension of $M$ to be even, say $\dim_\R M=2n$. The almost complex structure $J$ is called integrable when its Nijenhuis tensor, given by
	\begin{equation}\label{eq:nijenhuis}
		N_{J}(X,Y) = [X,Y]+J([JX,Y]+[X,JY])-[JX,JY],
	\end{equation}
	for $X,Y$ vector fields on $M$, is identically zero. An integrable almost complex structure is called simply a complex structure on $M$. According to the well-known Newlander-Nirenberg theorem \cite{N-N}, a complex structure on $M$ is equivalent to the existence of a holomorphic atlas on $M$, so that $(M,J)$ can be considered as a complex manifold of complex dimension $n$. 
	
	\indent Let $G$ be a connected Lie group with Lie algebra $\g$. An almost complex structure $J$ on $G$ is said to be left-invariant if left translations by elements of $G$ are complex maps. In this case $J$ is determined by the value at the identity of $G$. Thus, a left-invariant almost complex structure on $G$ amounts to an almost complex structure on its Lie algebra $\g$, that is, a real linear transformation $J$ of $\g$ satisfying $J^2 = -\Id$. Moreover, the former will be integrable if and only if the latter satisfies $N_J(x, y)=0$ for all $x, y$ in $\g$, with $N_J$ defined as in \eqref{eq:nijenhuis}.
	
	\indent Let $\g^\C$ denote the complexification of $\g$ and $\sigma$ the complex conjugation in $\g^\C$ with respect to the real form $\g$ of $\g^\C$. The $\C$-linear extension $J^\C\colon \g^\C\to\g^\C$ of any almost complex structure $J$ on $\g$ allows for a decomposition $\g^\C=\g^{1,0} \oplus \g^{0,1}$ of vector spaces, where $\g^{1,0}$ and  $\g^{0,1}$ are respectively the $i$-eigenspace and the $(-i)$-eigenspace  of $J^\C$, or alternatively  
	\begin{align*}
		\g^{1,0}=\{x-iJx : x\in \g\},\quad \g^{0,1} = \{x+iJx : x\in \g\}. 
	\end{align*}
	\noindent Note that these spaces satisfy $\g^{0,1} = \sigma(\g^{1,0})$. The integrability of $J$ is equivalent to $\g^{1,0}$ and $\g^{0,1}$ being both complex Lie subalgebras of $\g^\C$. Conversely, any complex Lie subalgebra $\m$ of $\g^\C$ satisfying $\g^\C=\m\oplus \sigma(\m)$ determines a complex structure $J$ on the real form $\g$:
	\begin{align*}
		Jx := 
		\begin{cases} 
			-  ix, & x \in \phantom{\sigma( }\m, \\
			\phantom{-} ix, & x\in\sigma(\m). 
		\end{cases} 
	\end{align*}    
	\noindent Such an endomophism $J$ commutes with $\sigma$ and hence defines an almost complex structure  $J\colon\g\to\g$, which is certainly integrable. A Lie subalgebra $\mathfrak m$ satisfying this condition is called an \textit{invariant complex subalgebra with respect to $\sigma$}. Thus we have that left-invariant complex structures on $G$, or complex structures on $\g$, are in correspondence with invariant complex subalgebras of $\g^\C$: every invariant complex subalgebra on $\g^\C$ can be identified with the $(-i)$-eigenspace of a complex structure $J$ on $\g$. This correspondence is of relevance in the context of Section \ref{section: SL(3,R)}.
	
	\indent Recall that two complex structures $J_1$ and $J_2$ on $\g$ are said to be equivalent if there exists a Lie algebra automorphism $\varphi \colon \g \to \g$ such that $J_1\circ \varphi=\varphi \circ J_2$. Note that we can restate this notion in the light of invariant complex subalgebras: if $\m_1$ and $\m_2$ are the invariant complex subalgebras corresponding respectively to the complex structures $J_1$ and $J_2$ then $J_1$ and $J_2$ are equivalent if and only if there exists a Lie algebra automorphism $\varphi$ of $\g^\C$ such that $\sigma\varphi=\varphi\sigma$ and $\varphi(\m_1)=\m_2$. With these definitions, the classification of complex structures on $\g$ is reduced to the classification of equivalence classes of invariant complex subalgebras.
	
	\medskip
	
	\subsection{Hypercomplex structures}\label{sect: hcx}
	
	\indent A hypercomplex structure on a smooth manifold $M$ is a triple of complex structures $\hcx$, $\alpha=1,2,3$, satisfying
	\begin{equation}\label{eq: quat} 
		J_1J_2=-J_2J_1=J_3. 
	\end{equation}
	Note in particular that $J_\alpha J_\beta =- J_\beta J_\alpha=J_\gamma$ for any cyclic permutation $(\alpha,\beta,\gamma)$ of $(1,2,3)$. It follows that such a manifold carries a 2-sphere of complex structures: they are
	\begin{equation*}
		J_a:=a_1J_1+a_2J_2+a_3J_3, \quad \text{where $a=(a_1,a_2,a_3)\in \mathbb{S}^2$}. 
	\end{equation*}
	\noindent Moreover, for any $p\in M$, the tangent space $T_p M$ has an $\mathbb{H}$-module structure and hence $\dim_\R M = 4n$ for some $n\in \N$. 
	
	\indent A hypercomplex structure $\hcx$ on $M$ determines a unique torsion-free connection $\nob$ such that $\nob J_\alpha = 0$ for $\alpha = 1, 2,3$, which is called  the Obata connection (see \cite{Ob}). As shown in \cite[equation (2.5)]{Sol}, an explicit expression for $\nob$ is given by:
	\begin{equation*}\label{eq:Obata}
		\nob_X Y=\tfrac12 \left([X,Y]+J_1[J_1 X,Y]-J_2[X,J_2 Y]+J_3 [J_1 X,J_2 Y]\right),
	\end{equation*}
	\noindent for $X$, $Y \in \X(M)$. 
	
	\indent Recall that the quaternionic general linear group $\mathrm{GL}(n, \mathbb{H})$ and its Lie algebra $\mathfrak{gl}(n, \mathbb{H})$ can be thought as subsets of $\mathrm{GL}(4n, \R)$ and $\mathfrak{gl}(4n, \R)$, respectively. To be explicit, if $\hcx$ is a hypercomplex structure on $\R^{4n}$ then
	\begin{gather}\label{eq: GL_nH} 
		\operatorname{GL}(n, \H) = \{T\in \operatorname{GL}(4n,\R) : TJ_\alpha =J_\alpha T \text{ for all } \alpha \}, \\
		\mathfrak{gl}(n,\mathbb{H})= \{T\in \mathfrak{gl}(4n,\R) : TJ_\alpha =J_\alpha T \text{ for all } \alpha \}. \notag
	\end{gather}
	
	After identifying any tangent space $T_pM$ of the hypercomplex manifold $M$ with $\R^{4n}$ via a choice of basis, it follows that the holonomy group of the Obata connection, $\operatorname{Hol}^{\mathrm{Ob}}:= \operatorname{Hol}(\nob)$, is contained in $\operatorname{GL}(n, \H)$ since each complex structure $J_\alpha$, $\alpha=1,2,3$, is $\nob$-parallel. 
	
	
	\indent As it is the case with complex structures, on Lie groups it makes sense to consider left-invariant hypercomplex structures, and they happen to be in correspondence with hypercomplex structures on Lie algebras. Naturally, the Obata connection corresponding to a left-invariant hypercomplex structure on a Lie groups is also left-invariant and it can be determined by its action on left-invariant vector fields, that is, on the Lie algebra. 
	
	We close this subsection recalling a description of left-invariant HKT structures on Lie groups. A left-invariant hyper-Hermitian structure on a Lie group $G$ is in correspondence with a pair $(\hcx, g)$ on the Lie algebra $\g$ of $G$, where $\hcx$ is a hypercomplex structure on $\g$ and $g$ is an inner product on $\g$ which is Hermitian with respect to $J_\alpha$ for all $\alpha$. It was proved in \cite[equation (3)]{DF} that the HKT condition (namely, the equality $\nabla^b_1=\nabla^b_2=\nabla^b_3$, where $\nabla^b_\alpha$ is the Bismut connection on $G$ associated to $(J_\alpha ,g)$) is equivalent to the following relation at the Lie algebra level:
	\begin{align} \label{eq: HKT}
		&g([J_1 x, J_1 y], z)+g([J_1 y,J_1 z],x)+g([J_1 z, J_1 x],y) \\ &=g([J_2 x, J_2 y], z)+g([J_2 y,J_2 z],x)+g([J_2 z, J_2 x],y),\nonumber \\ &=g([J_3 x, J_3 y], z)+g([J_3 y,J_3 z],x)+g([J_3 z, J_3 x],y)\nonumber 
	\end{align}
	\noindent for any $x$, $y$, $z\in \g$. The following result is an immediate consequence of the equations \eqref{eq: HKT}.
	
	\begin{lemma}\label{lem: HKT subalgebra}
		If $(\hcx, g)$ is an HKT structure on the Lie algebra $\g$ and $\h$ is a subalgebra of $\g$ invariant by the hypercomplex structure $\hcx$ then the restriction of $(\hcx, g)$ to $\h$ is again HKT.
	\end{lemma}
	
	\smallskip
	
	\subsection{Holonomy group of a linear connection} 
	
	As we are interested in studying the Obata holonomy of certain hypercomplex manifolds we collect here some well-known facts on holonomy groups and the Ambrose-Singer theorem that will be useful in subsequent sections.
	
	Let us consider a connected manifold $M$ equipped with a linear connection $\nabla$  and let us fix a point $p\in M$. If $\sigma:[0,1]\to M$ is a piecewise smooth loop based at $p$, the connection $\nabla$ gives rise to a parallel transport map $P_\sigma:T_pM \to T_pM$, which is linear and invertible. The holonomy group of $\nabla$ based at $p\in M$ is defined as 
	\[ \operatorname{Hol}_p(\nabla)=\{ P_\sigma\in \operatorname{GL}(T_pM) : \sigma \text{ is a loop based at } p\}.\] 
	It is well known that $\operatorname{Hol}_p(\nabla)$ is a Lie subgroup of $\operatorname{GL}(T_pM)$. 
	Since $M$ is connected, the holonomy groups based at two different points are conjugated, and therefore we can speak of the holonomy group of $\nabla$, denoted simply by $\operatorname{Hol}(\nabla)$. If $\dim M=n$, we can identify $\operatorname{Hol}(\nabla)$ with a Lie subgroup of $\operatorname{GL}(n,\R)$, after some choice of basis.
	The holonomy group need not be connected, and its identity component is denoted by $\operatorname{Hol}^0(\nabla)$; it is known as the restricted holonomy group of $\nabla$ and it consists of the parallel transport maps $P_\sigma$ where $\sigma$ is null-homotopic. Clearly, if $M$ is also simply connected then $\operatorname{Hol}(\nabla)=\operatorname{Hol}^0(\nabla)$ and it is connected.
	
	The Ambrose-Singer theorem provides a way to compute the holonomy group of a linear connection; indeed, it describes the Lie algebra $\mathfrak{hol}_p(\nabla)$ of $\operatorname{Hol}_p(\nabla)$ in terms of curvature endomorphisms $R_p(x,y)$ for $x,y\in T_pM$. 
	
	\begin{theorem}\cite{AmSi}
		The holonomy algebra $\mathfrak{hol}_p(\nabla)$ is the smallest subalgebra of $\mathfrak{gl}(T_pM)$ containing the endomorphims $P_\sigma^{-1}\circ R_p(x,y)\circ P_\sigma$, where $x,y$ run through $T_pM$, $\sigma$ runs through all piecewise smooth paths starting from $p$ and $P_\sigma$ denotes the parallel transport map along $\sigma$.
	\end{theorem}
	
	
	\smallskip
	
	Let us consider the particular case when $M=G$ is a Lie group with Lie algebra $\g$. A linear connection $\nabla$ on $G$ is said to be left-invariant if the left translations on $G$ are affine maps. As a consequence, if $X,Y$ are left-invariant vector fields then $\nabla_XY$ is also left-invariant. Therefore $\nabla$ is uniquely determined by a bilinear multiplication $\g \times \g \to \g$, still denoted by $\nabla$. We also denote by $\nabla_x\colon \g\to\g$ the endomorphism defined by left multiplication with $x\in\g$. In this case the Ambrose-Singer theorem takes the following form:
	
	\begin{theorem}\cite{Alek} \label{thm: hol}
		Let $\nabla$ be a left-invariant linear connection on the Lie group $G$, and let $\g$ be the Lie algebra of $G$. Then the holonomy algebra $\mathfrak{hol}(\nabla)$, based at the identity element $e\in G$, is the smallest Lie subalgebra of $\operatorname{End}(\g)$ containing the curvature endomorphisms $R(x,y)$ for any $x,y\in\g$, and closed under commutators with the left multiplication operators $\nabla_u\colon\g\to \g$.
	\end{theorem}
	
	We point out that the commutator $[\nabla_u,R(x,y)]$ can be interpreted as the covariant derivative $\nabla_u R(x,y)$, for $u,x,y\in\g$.
	
	\medskip
	
	\subsection{Construction of hypercomplex structures on complexified Lie algebras}
	
	We recall here a method developed in \cite{AS} to construct hypercomplex structures on the complexification of a Lie algebra equipped with a complex structure and a compatible decomposition.
	
	Let $\g$ be a real Lie algebra endowed with a complex structure $J$, and assume that there exist two Lie subalgebras $\g_+$ and $\g_-$ of $\g$ satisfying
	\begin{align*}
		\mathrm{ (1) } \; \g=\g_+\oplus \g_- \qquad \quad \text{ and } \quad \qquad  \mathrm{ (2) } \; J\g_+=\g_-.
	\end{align*}
	It follows that if $\dim \g=2n$ then $\dim\g_+=\dim \g_-=n$. In this case we say that $\g$ carries a complex product structure (see \cite{AS}). Another way to characterize such a structure is by defining an operator $E\in\operatorname{End} (\g)$ as: $Ex=x$ if $x\in\g_+$ and $Ex=-x$ if $x\in\g_-$. Since $E^2=\Id$ and the eigenspaces of $E$ corresponding to $\pm 1$, i.e. $\g_+$ and $\g_-$, are Lie subalgebras of $\g$, the endomorphism $E$ is called a \textit{product structure} on $\g$. Furthermore, it follows from $J\g_+=\g_-$ that $JE=-EJ$. Therefore, a complex product structure on $\g$ amounts to a pair $\{J,E\}$, where $J$ is a complex structure and $E$ is a product structure on $\g$ satisfying $JE=-EJ$. 
	
	Let us also recall from \cite{AS} the following important result.
	
	\begin{proposition} \cite[Section 5]{AS} \label{prop: existe ncp}
		Given a complex product structure $\{J,E\}$ on $\g$, there exists a unique torsion-free connection $\ncp$ on $\g$ such that $\ncp J=\ncp E=0$. The associated subalgebras $\g_+$ and $\g_-$ are totally geodesic and the induced torsion-free connections $\nabla^+$ on $\g_+$ and $\nabla^-$ on $\g_-$ are flat.
	\end{proposition}
	
	\begin{remark}
		It follows from Proposition \ref{prop: existe ncp} that $\g_+$ and $\g_-$ are not semisimple. Indeed, semisimple Lie algebras do not admit flat torsion-free connections, due to \cite[Theorem 7]{Chu}.
	\end{remark}
	
	The connection $\ncp$ is given explicitly by 
	\begin{align*}
		\ncp_xy=\tfrac{1}{4}\{&[x,y]-[Ex,Ey]+E[x,Ey]-E[Ex,y] -J[x,Jy]-J[Ex,Fy]+F[x,Fy]+F[Ex,Jy] \}
	\end{align*}
	for $x,\,y\in \g$, where $F:=JE$ is another product structure on $\g$. Another way to describe it is:
	\begin{equation}\label{eqs: ncp-casos}
		\begin{cases}
			\ncp_{x_+}y_+= -\pi_+J[x_+,Jy_+],\\
			\ncp_{x_-}y_-= -\pi_-J[x_-,Jy_-],\\
			\ncp_{x_+}y_-= \phantom{+} \pi_-[x_+,y_-],\\
			\ncp_{x_-}y_+= \phantom{+} \pi_+[x_-,y_+],
		\end{cases}
	\end{equation}
	for any $x_+,y_+\in\g_+,\;x_-,y_-\in\g_-$, where $\pi_{\pm}\colon \g\to\g_{\pm}$ are the projections. 
	Note that $\ncp_xy\in\g_{\pm}$ whenever $y\in\g_{\pm}$, for any $x\in\g$.
	
	\smallskip
	
	\indent Denote by $\hat{\g}:=(\g^\C)_\R$ the complexification of $\g$ when viewed as a real Lie algebra. As the following result shows, a complex product structure on $\g$ induces a hypercomplex structure on $\hat{\g}$.
	
	\begin{proposition} \label{prop: HCX} \cite[Theorem 3.3]{AS}
		If $\g$ admits a complex product structure $\{J,E\}$ then the Lie algebra $\hat{\g}$ carries a hypercomplex structure $\{J_1,J_2, J_3=J_1J_2\}$ given by
		\begin{gather}\label{eq: J1y2}
			J_1(x+iy) = J x+i J y, \quad x,y\in \g, \\
			J_2(x+iy) = \begin{cases} 
				\phantom{-}i(x+iy), & x,y\in\g_+, \\
				-i(x+iy),&  x,y\in\g_-.\nonumber
			\end{cases}
		\end{gather}
	\end{proposition}
	\noindent The associated Obata connection $\nob$ on $\hat{\g}$ is closely related to the torsion-free connection $\ncp$ on $\g$; indeed, it was shown in \cite{AS} that $\nob$ is the $\C$-linear extension of $\ncp$ to $\hat{\g}$, that is,
	\[ \nob_{x+iy}(z+iw)=(\ncp_x z-\ncp_y w)+i(\ncp_x w+\ncp_y z), \quad x,y,z,w\in \g.\]
	This implies that the Obata curvature tensor $R^{\operatorname{Ob}}$ is the $\C$-linear extension of $R^{\operatorname{CP}}$ and the same fact holds for the  covariant derivatives.
	The following result is an easy consequence of these facts and Theorem \ref{thm: hol}. We denote by $G$ (respectively, $\hat{G}$) the simply connected Lie group associated to $\g$ (respectively, $\hat{\g}$).
	
	\begin{proposition}\label{prop: complex-hol}
		The Obata holonomy algebra $\hol^\mathrm{Ob}$ of $\hat{G}$ is the complexification of the holonomy algebra $\hol^\mathrm{CP}$ of $G$. More precisely, $\hol^\mathrm{Ob}=((\hol^\mathrm{CP})^\C)_\R$.
	\end{proposition}
	
	\smallskip
	
	As mentioned before, hypercomplex manifolds with Obata holonomy contained in $\SL(n,\H)$ are special. Using Proposition \ref{prop: complex-hol} we will interpret this condition in terms of the holonomy of the connection $\ncp$. 
	
	Let $\g$ be a $2n$-dimensional Lie algebra equipped with a complex product structure $\{J,E\}$ and denote by $G$ the simply connected Lie group associated to $\g$. If $T\in \operatorname{Hol}^{\mathrm{CP}}$ then it follows from $\ncp J=0$ and $\ncp E=0$ that $T$ commutes with both $J$ and $E$. Choosing a basis $\{e_1,\ldots,e_n,Je_1,\ldots,Je_n\}$ of $\g$, where $e_j\in \g_+$ for all $j$, we obtain that $T=\left[\begin{array}{c|c} A & \boldsymbol{0} \cr \hline \boldsymbol{0} & A\cr \end{array}\right]$ for some $A\in \GL(n,\R)$. Hence, we arrive at the following result.  
	
	\begin{proposition}\label{prop: hol CP}
		The holonomy $\operatorname{Hol}^{\mathrm{CP}}$ of the left-invariant connection $\ncp$ on the $2n$-dimensional Lie group $G$ satisfies 
		\[\operatorname{Hol}^{\mathrm{CP}}\subseteq \GL(n,\R),\] where $\GL(n,\R)$ is considered as a subgroup of $\GL(2n,\R)$ via the diagonal inclusion
		\begin{equation}\label{eq: inclusion} 
			\GL(n,\R) \hookrightarrow \GL(2n,\R), \quad A\mapsto \left[\begin{array}{c|c} A & \boldsymbol{0} \cr \hline \boldsymbol{0} & A\cr \end{array}\right].  
		\end{equation}
	\end{proposition}
	\noindent As a consequence, the holonomy algebra $\hol^{\mathrm{CP}}$ is contained in $\gl(n,\R)$, considered as a Lie subalgebra of $\gl(2n,\R)$ via an inclusion similar to \eqref{eq: inclusion}.
	
	In the sequel, for short, we will denote $A^{\oplus 2}:=\left[\begin{smallarray}{c|c} A & \boldsymbol{0} \cr \hline \boldsymbol{0} & A\cr \end{smallarray}\right]$, and more generally, we will use $A\oplus B$ to denote the block-diagonal matrix $\left[\begin{smallmatrix}A&\\&B\end{smallmatrix}\right]$. 

Next we identify explicitly $\gl(n,\H)$ inside $\gl(4n,\R)$ as well as some special subalgebras such as $\gl(n,\C)$ and $\sl(n,\H)$. If we decompose $\g=\g_+\oplus \g_-$, let us consider a basis $\mathcal{B}_1=\{e_1,\ldots,e_n\}$ of $\g_+$. Then 
\begin{align*}
	\mathcal{B}_2=\{ie_1,\ldots,ie_n\}, \quad \mathcal{B}_3=\{Je_1,\ldots,Je_n\}, \quad\text{and} \quad\mathcal{B}_4=\{iJe_1,\ldots,iJe_n\}
\end{align*}
\noindent are $\R$-bases of $i\g_+$, $\g_-$, $i\g_-$ respectively, and $\mathcal{B}=\mathcal{B}_1\cup\mathcal{B}_2\cup\mathcal{B}_3\cup\mathcal{B}_4$ is a basis of $\hat{\g}$. With respect to this ordered basis, the matrices of the complex structures $J_1$ and $J_2$ from \eqref{eq: J1y2} are given by
\[ J_1=\left[\begin{array}{c|c|c|c}  & &-I_n & \cr \hline  & & & -I_n  \cr \hline  I_n & && \cr \hline & I_n && \cr\end{array}\right], \qquad 
J_2=\left[\begin{array}{c|c|c|c}  & -I_n& & \cr \hline  I_n& & &  \cr \hline   & && I_n\cr \hline &  & -I_n& \cr\end{array}\right],  \]
where $I_n$ denotes the $n\times n$ identity matrix. Considering $\{J_1,J_2,J_3:=J_1J_2\}$ as a hypercomplex structure on $\R^{4n}$, we have that 
\[ \gl(n,\H)=\{H(X,Y,Z,W)\in \gl(4n,\R): X,Y,Z,W\in \gl(n,\R)\}, \]
where
\begin{equation}\label{eq: H}  H(X,Y,Z,W)=\begin{bmatrix} 
		X & -Y & -Z & -W \\
		Y & \phantom{+} X & -W & \phantom{+}Z \\
		Z & \phantom{+} W & \phantom{+}X & -Y \\
		W & -Z & \phantom{+}Y & \phantom{+}X
	\end{bmatrix}.
\end{equation}
Furthermore, it is readily verified that 
\begin{equation}\label{eq: sl(n,H)}
	\sl(n,\H)=\{H(X,Y,Z,W)\in \gl(4n,\R): X\in \sl(n,\R), \, Y,Z,W\in \gl(n,\R)\}.
\end{equation} 

\begin{corollary} \label{coro: holonomies}
	\leavevmode
	\begin{enumerate} 
		\item[$\ri$] The Obata holonomy $\operatorname{Hol}^{\mathrm{Ob}}$ of the $4n$-dimensional simply connected Lie group $\hat{G}$ is contained in a proper subgroup of $\GL(n,\H)$, which is isomorphic to $\GL(n,\C)$. 
		\item[$\rii$] $\operatorname{Hol}^{\mathrm{Ob}}$ lies in $\SL(n,\H)$ if and only if the holonomy $\operatorname{Hol}^{\mathrm{CP}}$ of the $2n$-dimensional simply connected Lie group $G$ is contained in $\SL(n,\R)$, where  $\SL(n,\R)$ is considered as a Lie subgroup of $\GL(2n,\R)$ via \eqref{eq: inclusion}. In this case, $\operatorname{Hol}^{\mathrm{Ob}}$ is contained in a proper subgroup of $\SL(n,\H)$, isomorphic to $\SL(n,\C)$.
	\end{enumerate}
\end{corollary}

\begin{proof}
	Since $\hat{G}$ and $G$ are both simply connected, it will be enough to show the analogous statement for the corresponding holonomy algebras. 
	
	We keep the notation used previously. Let $T\in\hol^{\mathrm{Ob}}\subseteq \operatorname{End}(\hat{\g})$. It follows from $\hol^{\mathrm{CP}} \subseteq \gl(n,\R)$ and Proposition \ref{prop: complex-hol} that $T$ can be written in the basis $\mathcal{B}$ of $\hat{\g}$ as 
	\[ T=H(X,Y,0,0), \text{ where } X^{\oplus 2}, Y^{\oplus 2}\in \hol^{\operatorname{CP}}.\] 
	Note that $\mathfrak{v}:=\{H(X,Y,0,0): X,Y\in \gl(n,\R)\}$ is a Lie subalgebra of $\gl(n,\H)$, isomorphic to $\gl(n,\C)$, which is easily seen by mapping $H(X,Y,0,0)\in \gl(n,\H)$ to $X+iY\in\gl(n,\C)$. This proves $\ri$. 
	
	Assume now $\hol^{\mathrm{Ob}}\subseteq  \sl(n,\H)$ and let $X^{\oplus 2}\in \hol^{\mathrm{CP}}$. Then $T=H(X,0,0,0)\in\hol^{\mathrm{Ob}}$ and therefore it follows from \eqref{eq: sl(n,H)} that $X\in \sl(n,\R)$, thus $\hol^{\mathrm{CP}}\subseteq \sl(n,\R)$. Conversely, if  $\hol^{\mathrm{CP}} \subseteq \sl(n,\R)$ then it is clear that for any $T=H(X,Y,0,0)\in\hol^{\mathrm{Ob}}$ we have that $X,Y\in \sl(n,\R)$ so that $T\in \sl(n,\H)$. Moreover, $\hol^{\operatorname{Ob}}\subseteq \{H(X,Y,0,0): X,Y\in \sl(n,\R)\}\subseteq \mathfrak{v}$ and under the isomorphism between $\mathfrak{v}$ and $\gl(n,\C)$ we see that the subspace $\{H(X,Y,0,0): X,Y\in \sl(n,\R)\}$ is isomorphic to $\sl(n,\C)$. This proves $\rii$.
\end{proof}

\begin{corollary}
	The Obata holonomy of a simply connected $4n$-dimensional hypercomplex Lie group $\hat{G}$ whose Lie algebra $\hat{\g}$ carries a hypercomplex structure as in Proposition \ref{prop: HCX} is never equal to $\GL(n,\H)$ or $\SL(n,\H)$.
\end{corollary}

\begin{remark}
	In \cite[Example 6.3]{AT}, Proposition \ref{prop: HCX} is used to produce an 8-dimensional hypercomplex solvmanifold $(M, \{J_1, J_2, J_3\})$ with flat Obata connection and Obata holonomy group not contained in $\SL(n,\H)$ since the canonical bundle of the complex manifold $(M,J_2)$ is not holomorphically trivial. 
\end{remark}

\begin{remark}
	The holonomy $\operatorname{Hol}^{\mathrm{Ob}}$ of the Obata connection is reducible on $\hat{G}$. This follows from Proposition \ref{prop: complex-hol} and the fact that the holonomy $\operatorname{Hol}^{\mathrm{CP}}$ is reducible on $G$, as it is noted in Proposition \ref{prop: hol CP} that the Lie subalgebras $\g_{\pm}$ are $\operatorname{Hol}^{\mathrm{CP}}$-invariant. Thus, the holonomy group $\operatorname{Hol}^{\mathrm{Ob}}$ on $\hat{G}$ does not appear in the  Merkulov-Schwachh\"ofer classification of holonomy groups of irreducible torsion-free connections (see \cite{MeS}).
\end{remark}

\smallskip

\section{Non-existence of invariant hypercomplex structures on \texorpdfstring{$\SL(3,\R)$}{}} \label{section: SL(3,R)}

\indent In this section we show that there are no left-invariant hypercomplex structures on $\SL(3,\R)$: this is done in Subsection \ref{subsection: no hpcx}. To that aim, we use the classification of left-invariant complex structures $J$ on $\SL(3, \R)$ that was carried about in \cite{Sasaki}. We have found that this classification can be improved as it contains some flaws, so Subsection \ref{subsection: classification} is devoted to its revision. Subsection \ref{subsection: no hpcx} also contains a description of a \textit{non} left-invariant hypercomplex structure on $\SL(3,\R)$; see Remark \ref{obs: there is a hypercomplex structure}.  
\medskip

\subsection{Revision of the classification of left-invariant complex structures on \texorpdfstring{$\SL(3,\R)$}{}} \label{subsection: classification}

\indent We describe the classification of left-invariant complex structures on $\SL(3,\R)$ in terms of invariant complex subalgebras $\m$ of the complexification $\sl(3, \C)$ of $\sl(3, \R)$, as indicated in Section \ref{section: preliminares}. We need this classification as a intermediate step in our proof that $\SL(3,\R)$ has no left-invariant hypercomplex structures (see Theorem \ref{thm: SL(3,R) no tiene HCS} below). Let us fix such an $\m$ of $\g^{\C}$, and let
\begin{align} \label{eq: B y sigma B}
	\mathcal{B} = \{ U, X, Y, Z \} \quad \text{and} \quad \sigma(\mathcal{B}) = \{ U^{\sigma}, X^{\sigma}, Y^{\sigma}, Z^{\sigma} \}
\end{align}
\noindent denote a complex basis for $\m$ and the set of $\sigma$-conjugates of elements of $\mathcal{B}$, respectively. By means of specifying possible bases as in equation \eqref{eq: B y sigma B}, the following families of invariant complex subalgebras are described in \cite[equation (13)]{Sasaki}:
\begin{align} \label{eq: las ocho familias}
	\mathrm{I}_\lambda, \quad \mathrm{I}_\infty, \quad \mathrm{II}, \quad \mathrm{III}, \quad \mathrm{IV}_\mu, \quad\mathrm{V}, 
	\quad \mathrm{VII}_\lambda, \quad \mathrm{VII}_\infty \quad \mathrm{VIII}_\kappa, 
\end{align}
\noindent here $\lambda\in \C$ with $ |\lambda|\neq 1$, $\mu\in \C$, $\kappa\in \C$ with $\kappa\neq 1$.

In \cite[Theorem 1]{Sasaki} it is stated that every invariant complex subalgebra $\m$ is equivalent to one appearing in equation \eqref{eq: las ocho familias}, and that every one of these is equivalent to one and only one among
\begin{align} \label{eq: I, II, and III}
	\mathrm{I}_\lambda \text{ with $\lambda \in \C$ such that $| \lambda | < 1$}, \quad \mathrm{II}, \quad \mathrm{III};
\end{align}
\noindent in particular, no two subalgebras in equation \eqref{eq: I, II, and III} are equivalent. While we believe that this classification is correct in general lines, some intermediate steps appear to be wrong due to unintended miscalculations. In revising the argument we have found that, contrary to what was stated in \cite[Proposition 13]{Sasaki}, subalgebras $\mathrm{II}$ and $\mathrm{III}$ are actually equivalent (see Proposition \ref{prop: la 13 de sasaki pero corregida} below). Moreover, we believe that the classification up to equation \eqref{eq: las ocho familias}, proved in \cite[Propositions 4 through 8]{Sasaki}, is correct. The method used for establishing equivalence among subalgebras in equation \eqref{eq: las ocho familias} is computing the Lie brackets of the form $[\mathcal{B}, \mathcal{B}]$ and $[\mathcal{B}, \sigma(\mathcal{B})]$ for each subalgebra, where $\mathcal{B}$ and $\sigma(\mathcal{B})$ are as in equation \eqref{eq: B y sigma B}, and then using \cite[Lemma 1]{Sasaki} to establish that two subalgebras with the same Lie brackets must be equivalent. These Lie brackets are summarized in \cite[equations (14), (15), and (16)]{Sasaki}. However, we have found that some Lie brackets do not match the bases described in \cite[equation (13)]{Sasaki}. Fortunately, a slight modification of the bases of the subalgebras $\mathrm{IV}_\mu$ (with $\mu \neq 1$), $\mathrm{V}$, and $\mathrm{VII}_\infty$ is enough to ensure that the Lie bracket relations in \cite[equations (14), (15)]{Sasaki} hold. 
We describe them below in the notation used in \cite{Sasaki}: 
\begin{itemize}
	\item For the subalgebra $\mathrm{IV}_\mu$ (with $\mu\neq -1$), the following basis works:
	\begin{align*}
		U&=H_\alpha+2H_\beta+3e_\beta, \quad X=\frac{1}{\mu'}(H_\alpha+e_\alpha+\mu' e_{-\alpha}-e_\beta+e_\gamma), \\ 
		Y&= \overline{\mu'} e_\beta, \quad Z = \frac{\overline{\mu'}}{\mu'}(-e_\beta+e_\gamma).
	\end{align*}
	\noindent Here, $\mu' \in \C$ is defined such that $(\mu')^2 = 1 + \mu$. 
	\item For the subalgebra $\mathrm{V}$, the following basis works:
	\[  U=H_\alpha+2e_{-\alpha}+e_\beta, \quad X = \frac{1}{2} (e_\alpha - e_{-\alpha} - H_\alpha - e_\beta - e_\gamma), \quad
	Y=2 e_\beta, \quad Z = e_\beta + e_\gamma.\]  
	\noindent Note that this basis is very similar to the one in \cite[equation (13)]{Sasaki}, differing only in $U$. 
	\item For the subalgebra $\mathrm{VII}_\infty$, the following basis works:
	\begin{align*}
		U= H_\beta+e_{-\alpha}+2e_\beta, \quad X = H_\alpha+e_\alpha-e_{-\alpha}-e_\beta+e_\gamma, \quad
		Y= -e_\beta, \quad Z = e_\beta-e_\gamma.
	\end{align*}
\end{itemize}

\indent The rest of the classification goes on as in \cite[Propositions 9 through 12]{Sasaki}, with the exception of \cite[Proposition 13]{Sasaki}. Indeed, the subalgebras $\mathrm{II}$ and $\mathrm{III}$ \textit{are} equivalent. This can be seen by applying \cite[Lemma 1]{Sasaki} with the basis of $\mathrm{II}$ described in \cite[equation (13)]{Sasaki} and the following basis of $\mathrm{III}$:
\begin{align*} 
	U=\frac12 H_\alpha+H_\beta-\frac32 (e_\alpha+e_{-\alpha}-e_\beta+e_\gamma)&, \quad X=\frac12 H_\alpha+H_\beta+\frac12 (e_\alpha+e_{-\alpha}+e_\beta-e_\gamma), \\
	Y=-\frac12 (H_\alpha+e_\alpha-e_{-\alpha}-e_\beta+e_\gamma)&, \quad Z=\frac12(H_\alpha+e_\alpha-e_{-\alpha}+e_\beta-e_\gamma). \notag
\end{align*}
As this is a major contradiction from what was stated in \cite{Sasaki}, we collect this fact in the following result. 
\begin{proposition}\label{prop: la 13 de sasaki pero corregida}
	The invariant complex subalgebras $\mathrm{II}$ and $\mathrm{III}$ are equivalent.
\end{proposition}

\begin{remark} We point out that there is no basis of $\mathrm{III}$ satisfying \cite[equation (16)]{Sasaki}, due to Proposition \ref{prop: la 13 de sasaki pero corregida}.
\end{remark}

\indent The revised classification of invariant complex subalgebras of $\sl(3, \C)$ follows from \cite[Theorem 1]{Sasaki} and our corrections described above. Due to the correspondence between invariant complex subalgebras of $\sl(3, \C)$ and left-invariant complex structures on $\SL(3, \R)$, we have the following revision of Sasaki's classification.  
\begin{theorem}\label{thm: sasaki}
	Every left-invariant complex structure on $\SL(3,\R)$ is equivalent to one and only one of the complex structures corresponding with the invariant complex subalgebras $\mathrm{I}_\lambda$ $($with $\lambda \in \C$ such that $| \lambda | < 1)$ and $\mathrm{II}$.
\end{theorem} 
\indent 
We provide next more information about the Lie algebras $\mathrm{I}_\lambda$ and $\mathrm{II}$.

\begin{itemize}
	\item For the family $\mathrm{I}_\lambda$, the basis $\{U,X,Y,Z\}$ of the invariant complex subalgebra described in \cite[equation (13)]{Sasaki} is given by 
	\begin{gather*}
		U =\frac12 \matriz{1-\lambda & -i(\lambda+1) &0\\ i(\lambda+1)&1-\lambda&0\\0&0&-2(1-\lambda)},  \\
		X=\frac{1}{\sqrt{2}} \matriz{0&0&1\\0&0&i\\0&0&0}, \quad Y = \frac{1}{\sqrt{2}}\matriz{0&0&0\\0&0&0\\1&i&0}, \quad Z = \frac12 \matriz{1&i&0\\i&-1&0\\0&0&0}. \notag
	\end{gather*}  
	\item For the family $\mathrm{II}$, the basis $\{U,X,Y,Z\}$ of the invariant complex subalgebra described in \cite[equation (13)]{Sasaki} is given by 
	\begin{gather*}
		U=\frac12 \matriz{-1&-3i&0\\3i&-1&0\\0&0&2}, \quad X=\frac{1}{\sqrt{2}}\matriz{0&0&1\\0&0&i\\1&-i&0}, \quad
		Y=\frac{1}{\sqrt{2}}\matriz{0&0&0\\0&0&0\\1&i&0}, \quad Z=\frac12 \matriz{1&i&0\\i&-1&0\\0&0&0}.
	\end{gather*}  
\end{itemize}

Lie brackets of the form $[\mathcal{B}, \mathcal{B}]$ and $[\mathcal{B}, \sigma(\mathcal{B})]$ for the bases described above are listed in Tables \ref{tab: lie brackets for I and II} and \ref{tab: lie brackets for I and II sigma}. These relations will be used in the proof of Theorem \ref{thm: SL(3,R) no tiene HCS}. 
\begin{table}[H]
	\centering
	\renewcommand{\arraystretch}{1.5}
	\begin{tabular}{|c|cccccc|} 
		\hline
		& $[U,X]$ & $[U,Y]$ & $[U,Z]$ & $[X,Y]$ & $[X, Z]$ & $[Y, Z]$ \\
		\hline 
		\multirow[t]{4}{*}{$\mathrm{I}_\lambda$} & $(2-\lambda)X$& $(2\lambda-1)Y$ & $(\lambda+1)Z$ & $Z$ & $0$ & $0$ \\    
		\hline 
		\multirow[t]{4}{*}{$\mathrm{II}$} & $0$ & $3 Y$ & $3 Z$ & $Z$ & $Y$ & $0$ \\    
		\hline 
	\end{tabular}
	\caption{Lie brackets of the form $[\mathcal{B}, \mathcal{B}]$ for the invariant complex subalgebras $\mathrm{I}_\lambda$ and $\mathrm{II}$.}
	\label{tab: lie brackets for I and II}
\end{table}
\begin{table}[H]
	\centering
	\renewcommand{\arraystretch}{1.5}
	\begin{tabular}{|c|cccccccc|} 
		\hline
		& $[U,X^{\sigma}]$ & $[U,Y^{\sigma}]$ & $[U,Z^{\sigma}]$ & $[X,X^{\sigma}]$ & $[X,Y^{\sigma}]$ & $[X, Z^{\sigma}]$ & $[Y, Z^{\sigma}]$ & $[Z, Z^{\sigma}]$ \\
		\hline 
		\multirow[t]{4}{*}{$\mathrm{I}_\lambda$} & $(1-2\lambda)X^{\sigma}$& $(\lambda-2)Y^{\sigma}$ & $-(\lambda+1)Z^{\sigma}$ & $0$ & $\frac{U+\lambda U^\sigma}{1-|\lambda|^2}$ & $-X^\sigma$ & $Y^\sigma$ & $\frac{(1-\overline{\lambda})U-(1-\lambda)U^\sigma}{1-|\lambda|^2}$ \\ 
		\hline 
		\multirow[t]{4}{*}{$\mathrm{II}$} & $6Y - 3 X^{\sigma}$ & $0$ & $- 3 Z^{\sigma}$ & $Z-Z^\sigma$ & $-\frac{U+2U^\sigma}{3}$ & $Y-X^\sigma$ & $Y^\sigma$ & $\frac{U-U^\sigma}{3}$ \\ 
		\hline 
	\end{tabular}
	\caption{Lie brackets  of the form $[\mathcal{B}, \sigma(\mathcal{B})]$ for the invariant complex subalgebras $\mathrm{I}_\lambda$ and $\mathrm{II}$.\\
		For both subalgebras we have that $[U, U^{\sigma}] = [Y, Y^{\sigma}] = 0$. }
	\label{tab: lie brackets for I and II sigma}
\end{table}

\indent For the sake of completeness we describe below the left-invariant complex structures corresponding to the invariant complex subalgebras $\mathrm{I}_\lambda$ and $\mathrm{II}$. We achieve this by describing their matrix form with respect to the ordered basis of $\sl(3,\R)$ given by
\begin{align*}
	\{E_{11}-E_{33},E_{22}-E_{33}, E_{21}, E_{12}, E_{31}, E_{32}, E_{13}, E_{23}\}.
\end{align*}
\noindent Here, $E_{ij}$ denotes the matrix with $1$ in the entry $(i,j)$ and $0$ in the others. Decomposing each element into its real and imaginary part and writing $\lambda=a+ib$, it is straightforward to verify that the corresponding left-invariant complex structures $J_{\mathrm{I}_\lambda}$ and $J_{\mathrm{II}}$ are 
\renewcommand{\arraystretch}{1.5} 
\begin{equation*}
	J_{\mathrm{I}_\lambda}=\frac{1}{a^2+b^2-1}\left[\begin{array}{cccc} 
		b&b&1-a&-a^2-b^2+a\\ b&b& a^2+b^2-a&-1+a\\ -a-1& -a^2-b^2-a&-b&\phantom{+}b\\ a^2+b^2+a&a+1&\phantom{+}b&-b
	\end{array}\right]\oplus \matriz{&0&-1&0&0\\&1&0&0&0\\&0&0&0&-1\\&0&0&1&0},
\end{equation*}  
\renewcommand{\arraystretch}{1.5}
\[J_{\mathrm{II}}=\left[\begin{array}{cccc}0&0&-\frac13&-\frac23\\0&0&\frac23&\frac13\\-1&-2&0&0\\2&1&0&0\end{array}\right] \oplus \matriz{ 0&-1&0&-2\\1&0&-2&0\\0&0&0&-1\\0&0&1&0}.\]
\renewcommand{\arraystretch}{1}

\smallskip

\subsection{Non-existence of left-invariant hypercomplex structures on \texorpdfstring{$\SL(3,\R)$}{}} \label{subsection: no hpcx}

\indent As an application of Theorem \ref{thm: sasaki} we find that $\SL(3, \R)$ does not admit left-invariant hypercomplex structures.
\begin{theorem} \label{thm: SL(3,R) no tiene HCS}
	The Lie group $\SL(3,\R)$ does not admit any left-invariant hypercomplex structure.
\end{theorem}
\begin{proof}
	Left-invariant hypercomplex structures on $\SL(3, \R)$ are in correspondence with hypercomplex structures on $\sl(3, \R)$, so we work entirely on the Lie algebra level as usual. We prove the result by means of contradiction. It is straightforward to verify that a hypercomplex structure on a Lie algebra $\g$ is equivalent to a pair $(\mathfrak m, J)$, where $\mathfrak m$ is an invariant complex subalgebra of $\g^\C$ and $J$ is a $\C$-linear endomorphism of $\g^\C$ satisfying
	\begin{equation}\label{eq: equiv hcpx subalgebra} 
		J^2=-\Id,\quad N_J=0, \quad \sigma J=J\sigma, \quad \text{and} \quad J(\mathfrak m)=\sigma(\mathfrak m).
	\end{equation} Thus, we fix such an invariant complex subalgebra $\m$ of $\sl(3,\C)$ and we show that there is no $\C$-linear isomorphism $J\colon \sl(3, \C) \to \sl(3, \C)$ such that \eqref{eq: equiv hcpx subalgebra} holds.
	We achieve this by performing a case-by-case analysis based on Theorem \ref{thm: sasaki}, depending on whether the invariant complex subalgebra $\m$ is equivalent to $\mathrm{I}_{\lambda}$ or $\mathrm{II}$. Notice that such a $J$ must be of the form 
	\begin{align*}
		J = \matriz{
			0&0&0&0&\overline{\lambda_1}&\overline{\lambda_5}&\overline{\lambda_9}&\overline{\lambda_{13}}\\ 
			0&0&0&0&\overline{\lambda_2}&\overline{\lambda_6}&\overline{\lambda_{10}}&\overline{\lambda_{14}}\\
			0&0&0&0&\overline{\lambda_3}&\overline{\lambda_7}&\overline{\lambda_{11}}&\overline{\lambda_{15}}\\
			0&0&0&0&\overline{\lambda_4}&\overline{\lambda_8}&\overline{\lambda_{12}}&\overline{\lambda_{16}}\\
			\lambda_1&\lambda_5&\lambda_9&\lambda_{13}&0&0&0&0\\
			\lambda_2&\lambda_6&\lambda_{10}&\lambda_{14}&0&0&0&0\\
			\lambda_3&\lambda_7&\lambda_{11}&\lambda_{15}&0&0&0&0\\
			\lambda_4&\lambda_8&\lambda_{12}&\lambda_{16}&0&0&0&0},
	\end{align*}
	\noindent for some $\lambda_j \in \C$, $1 \leq j \leq 16$, with respect to a complex basis
	\begin{align} \label{eq: base de sl(3,c)}
		\{U, X, Y, Z, U^{\sigma}, X^{\sigma}, Y^{\sigma}, Z^{\sigma} \}
	\end{align}
	\noindent  of $\sl(3, \C)$ such that $\{U, X, Y, Z \}$ is a complex basis of $\m$ as in Table \ref{tab: lie brackets for I and II}, and it must verify that $C := J^2 + \Id$ is identically zero. Note that the condition $C \equiv 0$ is independent of the invariant complex subalgebra $\m$. A useful observation which will be used repeatedly in the proof is that
	\begin{equation} \label{eq: c11}
		|\lambda_1|^2 + |\lambda_5|^2 + |\lambda_9|^2 + |\lambda_{13}|^2 \neq 0,
	\end{equation}
	\noindent i.e., the set $\{\lambda_l : l = 1, 5, 9, 13 \}$ has a nonzero element, which follows from the fact that    
	\begin{align*} 
		C_{11} = |\lambda_1|^2 + \overline{\lambda_5}\lambda_2 + \overline{\lambda_9}\lambda_3 + \overline{\lambda_{13}}\lambda_4 + 1
	\end{align*}
	\noindent must be zero. By notational convenience, we rename the basis in equation \eqref{eq: base de sl(3,c)} as $\{e_1, \ldots, e_8\}$. In the same vein, if $N := N_J$, let us call 
	\begin{align*}
		N_{ijk} := \la N(e_i, e_j), e_k\ra \quad \text{for all $1 \leq i, j , k \leq 8$}. 
	\end{align*}
	\noindent For $J$ to be integrable, it must be the case that $N_{ijk}$ be zero for all $1 \leq i,j,k \leq 8$. We now split the proof in two cases\footnote{We have double-checked all the following steps with symbolic manipulation software, as computations  are often very challenging.}.  \\ 
	
	\indent $\ri$ \underline{Case $\m = \mathrm{I}_{\lambda}$}. We begin by noticing that
	\begin{align} \label{eq: c66}
		C_{66} = \overline{\lambda_5} \lambda_2 + | \lambda_6 |^2 + \overline{\lambda_7} \lambda_{10} + \overline{\lambda_8} \lambda_{14} + 1;
	\end{align}    
	\noindent this expression must be equal to zero. Note that
	\begin{gather*}
		N_{184} = -2 \lambda_1 \overline{\lambda_{16}} (\overline{\lambda} + 1), \quad N_{284} = -2 \lambda_5    \overline{\lambda_{16}} (\overline{\lambda} + 1), \\ 
		N_{384} = -2 \lambda_9 \overline{\lambda_{16}} (\overline{\lambda} + 1), \quad N_{484} = -2 \lambda_{13} \overline{\lambda_{16}} (\overline{\lambda} + 1).
	\end{gather*}
	\noindent It follows from equation \eqref{eq: c11} and the condition $| \lambda | < 1$ that $\lambda_{16} = 0$. We assume from now on that we have set $\lambda_{16}=0$ in $J$, and thus in all $N_{ijk}$'s. Now, 
	\begin{gather*} 
		N_{183} = -3 \lambda_1 \overline{\lambda_{15}}, \quad N_{283} = -3 \lambda_5 \overline{\lambda_{15}}, \quad N_{383} = -3 \lambda_9 \overline{\lambda_{15}}, \quad N_{483} = -3 \lambda_{13} \overline{\lambda_{15}}, \\
		N_{164} = -3\lambda_1 \overline{\lambda_8}, \quad N_{264} = -3 \lambda_5 \overline{\lambda_8}, \quad N_{364} = -3 \lambda_9 \overline{\lambda_8}, \quad N_{464} = -3\lambda_{13} \overline{\lambda_8}.
	\end{gather*}
	\noindent It follows from equation \eqref{eq: c11} again that $\lambda_8 = \lambda_{15} = 0$. Once more, we assume from now on that we have set $\lambda_8 = \lambda_{15} = 0$ in $J$, and thus in all $N_{ijk}$'s. Then,
	\begin{gather*}
		N_{163} = 2\lambda_1 \overline{\lambda_7} (\overline{\lambda} - 2), \quad N_{263} = 2\lambda_5 \overline{\lambda_7} (\overline{\lambda} - 2),\\
		N_{363} = 2\lambda_9 \overline{\lambda_7} (\overline{\lambda}-2), \quad N_{463} = 2\lambda_{13} \overline{\lambda_7} (\overline{\lambda} - 2).
	\end{gather*}
	\noindent Again, it follows from equation \eqref{eq: c11} and the condition $| \lambda | < 1$ that $\lambda_7 = 0$. As before, we assume from now on that we have set $\lambda_7 = 0$ in $J$, and thus in all $N_{ijk}$'s. Finally, observe that
	\begin{align*}
		N_{461} = \overline{\lambda_5} \lambda_{13} (\overline{\lambda} - 2), \quad N_{462} = -\lambda_{13} \overline{\lambda_6} (\overline{\lambda}+1)-1. 
	\end{align*}
	\noindent Once more, it follows from equation \eqref{eq: c11} and the condition $| \lambda | < 1$ that $\lambda_5 = 0$ and $\lambda_6 \neq 0$. Summing up, we have found out that $\lambda_5 = \lambda_7 = \lambda_8 = 0$ and $\lambda_6 \neq 0$, so equation \eqref{eq: c66} reads 
	\begin{align*}
		C_{66} = | \lambda_6 |^2 + 1 > 0.
	\end{align*}
	\noindent This is the contradiction we were looking for. Therefore there is no hypercomplex structure on $\sl(3,\R)$ such that $\m = \mathrm{I}_\lambda$. \\
	
	\indent $\rii$ \underline{Case $\m = \mathrm{II}$}. We begin by computing
	\begin{gather*}
		N_{172} - N_{184} = 6 \lambda_1 (\overline{\lambda_{16}} - \overline{\lambda_{10}}), \quad N_{272} - N_{284} = 6 \lambda_5    (\overline{\lambda_{16}} - \overline{\lambda_{10}}), \\
		N_{372} - N_{384} = 6 \lambda_9 (\overline{\lambda_{16}} - \overline{\lambda_{10}}), \quad N_{472} - N_{484} = 6 \lambda_{13} (\overline{\lambda_{16}} - \overline{\lambda_{10}}), \\
		N_{174} - N_{182} = 6 \lambda_1 (\overline{\lambda_{14}} - \overline{\lambda_{12}}), \quad N_{274} - N_{282} = 6 \lambda_5    (\overline{\lambda_{14}} - \overline{\lambda_{12}}), \\
		N_{374} - N_{382} = 6 \lambda_9 (\overline{\lambda_{14}} - \overline{\lambda_{12}}), \quad N_{474} - N_{482} = 6 \lambda_{13} (\overline{\lambda_{14}} - \overline{\lambda_{12}}).
	\end{gather*}
	\noindent All these expressions must vanish. Then, it follows from equation \eqref{eq: c11} that $\lambda_{16} = \lambda_{10}$ and $\lambda_{14} = \lambda_{12}$. We assume from now on that we have set $\lambda_{16} = \lambda_{10}$ and $\lambda_{14} = \lambda_{12}$ in $J$, and thus in all $N_{ijk}$'s. We get from the vanishing of
	\begin{align*}
		N_{485} = \frac{1}{3} (3 \lambda_9 + 2\lambda_{15}) \overline{\lambda_{12}} - \frac13 \overline{\lambda_{15}} \lambda_{12} + \frac13 | \lambda_{10} |^2 + 3| \lambda_{13} |^2 + \frac{1}{3}
	\end{align*}
	\noindent that $\lambda_{12}\neq 0$, and from this and the vanishing of
	\begin{align*}
		N_{488} = 2 |\lambda_{12}|^2 + 6 \overline{\lambda_{13}} \lambda_{10}
	\end{align*}
	\noindent we obtain that $\lambda_{10} \neq 0$ as well. Now, the vanishing of 
	\begin{align*}
		N_{136} = \frac23 \lambda_2^2-6\lambda_{10}^2
	\end{align*}
	\noindent yields $\lambda_2 = 3 \varepsilon \lambda_{10}$ for some $\varepsilon = \pm 1$, the vanishing of
	\begin{align*}
		N_{342} = -\frac{2}{3} \overline{\lambda_2} \lambda_{12} - 2 | \lambda_{10} |^2 = - 2 \varepsilon \overline{\lambda_{10}} \lambda_{12} - 2| \lambda_{10}|^2 
	\end{align*}
	\noindent yields $\lambda_{12} = - \varepsilon \lambda_{10}$. Finally, the vanishing of
	\begin{align*}
		N_{344} = - \frac{2}{3} \lambda_{12} \overline{\lambda_4} - 2 \overline{\lambda_{12}} \lambda_{10} = \frac{2}{3} \varepsilon \overline{\lambda_4} \lambda_{10} + 2 \varepsilon | \lambda_{10} |^2
	\end{align*}
	\noindent together with the condition $\lambda_{10} \neq 0$ yields $\lambda_4 = -3 \lambda_{10}$, and the vanishing of
	\begin{gather*}
		N_{134} = \frac{1}{3} \overline{\lambda_4} \lambda_2 + 3 \overline{\lambda_8} \lambda_{10} + \overline{\lambda_{12}} \lambda_4 + 3 \overline{\lambda_{10}} \lambda_{12} = 3\left( - \varepsilon | \lambda_{10} |^2 + \overline{\lambda_8} \lambda_{10} \right), \\
		N_{144} = - \frac{1}{3} | \lambda_4|^2 + 3 \overline{\lambda_8} \lambda_{12} - \overline{\lambda_{12}} \lambda_2 + 3 |\lambda_{10} |^2 - 3 = 3(|\lambda_{10}|^2 - \varepsilon \overline{\lambda_{8}} \lambda_{10} - 1)
	\end{gather*} 
	\noindent yields the contradiction $0 = N_{134} \varepsilon + N_{144} = - 3$. Therefore there is no hypercomplex structure on $\sl(3,\R)$ such that $\m = \mathrm{II}$ and this concludes the proof. 
\end{proof}

\begin{remark} \label{obs: there is a hypercomplex structure}
	Although there is no \textit{left-invariant} hypercomplex structure, one can produce a hypercomplex structure on $\SL(3,\R)$. The argument is valid for $\SL(2n+1, \R)$ for arbitrary $n \in \N$ and it follows from the Iwasawa decomposition of this Lie group.  Recall that the Iwasawa decomposition of a non-compact semisimple Lie group $G$ ensures that it is diffeomorphic to a product $K \times A \times N$ of a compact Lie subgroup $K$, an abelian Lie subgroup $A$, and a nilpotent Lie subgroup $N$ of $G$; moreover, for $G = \SL(2n+1,\R)$, $K$ is isomorphic to $\mathrm{SO}(2n+1)$, and the product $A \times N$ is simply connected and has dimension\footnote{It is well-known that $A$ is the set of diagonal matrices with positive entries and determinant equal to $1$, and $N$ is the unipotent group consisting of upper triangular matrices with $1$'s on the diagonal.} $t :=n (2n+3)$  and thus is diffeomorphic to $\R^t$. It follows from \cite[Theorem 4.2]{Joyce} and \cite[Section 4.5]{DT1} that $\mathrm{SO}(2n+1)\times T^n$ admits a left-invariant hypercomplex structure; moreover, also the non-compact group $\mathrm{SO}(2n+1)\times \R^n$ admits one (since they have the same Lie algebra). As $t>n$, we see that $\SL(2n+1, \R)$ is diffeomorphic to the product
	\begin{align*}
		\underbrace{\mathrm{SO}(2n+1) \times \R^n}_{:= M} \times \; \mathbb{R}^{t-n}.
	\end{align*}
	\noindent Note that both $\dim M$ and $t-n$ are equal to $2n^2+2n$, which is a multiple of $4$ for all $n$. Since both $M$ and $\R^{t-n}$ admit hypercomplex structures, it follows that $\SL(2n+1,\R)$ admits a hypercomplex structure, arising as the product of the corresponding structures on each factor. 
	
	Furthermore, since $\mathrm{SO}(2n+1)\times \R^n$ admits an HKT metric (see for instance \cite[Section 2.3]{GP}) and $\R^{t-n}$ admits an HKT metric (actually, hyper-Kähler), the product hypercomplex structure with the product metric is an HKT structure on $\mathrm{SL}(2n+1,\R)$, which is not left-invariant.
\end{remark}

\begin{remark} 
	We point out that, according to \cite[Theorem C]{Bo}, the Lie group $\SL(2n+1,\R)$ has uniform lattices, that is, discrete subgroups $\Gamma$ such that the quotient $M_\Gamma:=\Gamma\backslash \SL(2n+1,\R)$ is compact.
	However, since the hypercomplex structure on $\SL(2n+1,\R)$ constructed in Remark \ref{obs: there is a hypercomplex structure} is not left-invariant, it does not necessarily descend to $M_\Gamma$. 
\end{remark}

\smallskip


\section{Existence of a left-invariant hypercomplex structure on \texorpdfstring{$\SL(2n+1,\C)$}{}} \label{section: SL(2n+1,C)}

In this section we describe a complex product structure on $\SL(2n+1,\R)$ and then apply Proposition \ref{prop: HCX} to construct a hypercomplex structure on $\SL(2n+1,\C)$ from it (see Corollary \ref{cor: hpcx underlying}). No restrictions are placed on $n \in \N$. We point out that a similar result was obtained in \cite{IT} using a different approach. At the end of the section we prove that there is no left-invariant HKT metric compatible with this hypercomplex structure (see Proposition \ref{prop: there is no compatible HKT metric}).  

\indent For any $m\in \N$, the canonical inclusion of the general linear group $\GL(m,\R)$ into $\SL(m+1,\R)$ given by
\begin{align*}
	\GL(m,\R) \ni X \mapsto \left[\begin{array}{c|c} X & \boldsymbol{0} \cr \hline \boldsymbol{0} & (\det X)^{-1}\cr \end{array}\right]\in \SL(m+1,\R),
\end{align*}
\noindent is a Lie group monomorphism and induces naturally an inclusion of their Lie algebras as 
\begin{align} \label{eq: inclusion of gl in sl}
	\gl(m,\R) \ni X \mapsto \left[\begin{array}{c|c} X & \boldsymbol{0} \cr \hline \boldsymbol{0} & -\tr X \cr  \end{array}\right]\in \sl(m+1,\R).
\end{align}
\noindent Setting
\begin{align*}
	\a:=\left\{\left[\begin{array}{c|c} \boldsymbol{0} & \boldsymbol{0}
		\cr \hline v^t & 0\cr\end{array}\right]:\;v\in\R^m\right\} \quad \text{and}\quad 
	\b:=\left\{\left[\begin{array}{c|c} \boldsymbol{0} & w\cr \hline \boldsymbol{0} & 0\cr\end{array}\right]:\;w\in\R^m\right\},
\end{align*}
\noindent and taking $\gl(m,\R)$ as the image of the inclusion above, we have 
the vector space decomposition 
\begin{align} \label{eq: decomposition of sl in gl a b}
	\sl(m+1,\R)=\gl(m,\R)\oplus\a\oplus\b.
\end{align}
\noindent The Lie brackets in $\sl(m+1,\R)$ can be described easily if one makes the identification
\begin{align*}
	\sl(m+1,\R) \ni X = \left[\begin{array}{c|c} A & w\cr \hline v^t & -\tr A\cr \end{array}\right] \quad \longleftrightarrow \quad (A,v,w) \in \gl(m,\R)\oplus\a\oplus\b; 
\end{align*}
\noindent they turn out to be the following:  
\begin{align*}
	[(A,0,0),(A',0,0)]& =([A,A'],0,0), \\
	[(A,0,0),(0,v',0)]& =(0,-A^tv'-(\tr A)v',0), \\
	[(A,0,0),(0,0,w')]& =(0,0,Aw'+(\tr A)w'), \\
	[(0,v,0),(0,v',0)]& =(0,0,0), \\
	[(0,v,0),(0,0,w')]& =(-w'\cdot v^t,0,0),\\
	[(0,0,w),(0,0,w')]& =(0,0,0).
\end{align*} 
\noindent Note that ``$w' \cdot v^t$'' is the outer product of the column-vector $w'$ with the row-vector $v^t$. In particular, both $\a$ and $\b$ are abelian Lie subalgebras of $\sl(m+1,\R)$, isomorphic to $\R^m$. From now on, we will consider the case $m=2n$ for $n\geq 1$.

It has been proved in \cite[Proposition 3.5]{AS} that $\gl(2n,\R)$ carries a natural complex product structure $\{J',E'\}$, where the complex structure $J'$ and the product structure $E'$ are given by right multiplication with the matrices
\begin{align*}
	J_0=\left[\begin{array}{c|c} \boldsymbol{0} & -I_n\cr \hline I_n & \boldsymbol{0}\cr\end{array}\right] \quad \text{and} \quad E_0=\left[\begin{array}{c|c}I_n & \boldsymbol{0}\cr  \hline\boldsymbol{0} & -I_n\cr\end{array}\right],
\end{align*}
\noindent respectively. We now produce a complex product structure $\{J, E\}$ on
$\sl(2n+1,\R)$ by extending this complex product structure on the Lie subalgebra $\gl(2n,\R)$ (see Proposition \ref{prop: complex product structure} below). We achieve this by defining the endomorphisms $J$ and $E$ of $\sl(2n+1,\R)$ as
\begin{align*}
	J(A,v,w) := (AJ_0,-J_0v,-J_0w), \quad E(A,v,w) := (AE_0,E_0v,-E_0w). 
\end{align*} 
\noindent Certainly $J^2 = -\Id$ and $E^2 = \Id$, so that $J$ and $E$ are, respectively, an almost complex structure and an almost product structure on $\sl(2n+1,\R)$. It is clear that $J$ and $E$ anticommute. Note that, in terms of matrices, $J$ and $E$ act as follows:
\begin{gather*}
	J\left[\begin{array}{cc|c} A & C & w_1 \cr B & D & w_2 \cr \hline v_1^t & v_2^t & -\tr(A+D)\cr \end{array}\right]=
	\left[\begin{array}{cc|c} C & -A & w_2 \cr D & -B & -w_1 \cr \hline v_2^t & -v_1^t & \tr(B-C) \cr \end{array}\right],  \\
	E\left[\begin{array}{cc|c} A & C & w_1 \cr B & D & w_2 \cr \hline v_1^t & v_2^t & -\tr(A+D)\cr \end{array}\right]= 
	\left[\begin{array}{cc|c} A & -C & -w_1 \cr B & -D & w_2 \cr  \hline v_1^t & -v_2^t & \tr(D-A)\cr \end{array}\right].
\end{gather*}  
\noindent In particular, the eigenspaces $\s_{\pm}$ with eigenvalues $\pm 1$ of $E$  are given by:
\begin{gather*}
	\s_+=\left\{\left[\begin{array}{cc|c} A & \boldsymbol{0} & \boldsymbol{0}\cr B & \boldsymbol{0} & w_2 \cr \hline  v_1^t & \boldsymbol{0} & -\tr A\cr \end{array}\right] :  A,B\in \gl(n,\R), v_1,w_2\in\R^n\right\}, \\
	\s_-=\left\{\left[\begin{array}{cc|c}\boldsymbol{0} & C & w_1\cr \boldsymbol{0} & D & \boldsymbol{0}\cr \hline \boldsymbol{0} & v_2^t & -\tr D\cr\end{array}\right] : C,D\in \gl(n,\R), v_2,w_1\in\R^n\right\}. 
\end{gather*} 
\noindent Note that $\s_+$ and $\s_-$ are both subalgebras of $\sl(2n+1,\R)$, as can be readily verified, from where it follows that $E$ is indeed a product structure on $\sl(2n+1,\R)$. Moreover, $J$ is integrable; in fact, we will show that the relation 
\begin{equation}\label{eq: J es integ} 
	J[X,Y]=[JX,Y]+[X,JY]+J[JX,JY]
\end{equation}
\noindent holds for all $X,Y\in\sl(2n+1,\R)$, and this is equivalent to $N_J=0$. We verify \eqref{eq: J es integ} by considering several cases: 

\medskip

\indent $(1)$ Take $X=(A,0,0) \in \gl(2n,\R)$ and $Y=(A',0,0) \in \gl(2n,\R)$. Then \eqref{eq: J es integ} holds because $J'$ is integrable on $\gl(2n,\R)$. \\
\indent $(2)$ Take $X=(A,0,0)\in\gl(2n,\R)$ and $Y=(0,v',0)\in\a$. Then  \eqref{eq: J es integ} holds because 
\begin{align*}
	J[(A,0,0),(0,v',0)] &= J(0,-A^tv'-(\tr A)v',0)=(0, J_0A^tv'+(\tr A)J_0v',0), \\
	[J(A,0,0),(0,v',0)] &= [(AJ_0,0,0),(0,v',0)]=(0,J_0A^tv'-(\tr AJ_0)v',0), \\
	[(A,0,0),J(0,v',0)] &= [(A,0,0),(0,-J_0v',0)]=(0,A^tJ_0v'+(\tr A)J_0v',0), 
\end{align*} 
\noindent and, as $J_0^t = -J_0$, 
\begin{align*}
	J[J(A,0,0),J(0,v',0)] &= J[(AJ_0,0,0),(0,-J_0v',0)]\\
	&= J(0,-J_0A^tJ_0v'+(\tr AJ_0)J_0v',0)\\
	&= (0,-A^tJ_0v'+(\tr AJ_0)v',0).
\end{align*}
\indent $(3)$ Take $X=(A,0,0)\in\gl(2n,\R)$ and $Y=(0,0,w')\in\b$. Then \eqref{eq: J es integ} holds because
\begin{align*}
	J[(A,0,0),(0,0,w')] &= J(0,0,Aw'+(\tr A)w')=(0,0,-J_0Aw'-(\tr A)J_0w'), \\
	[J(A,0,0),(0,0,w')] &= [(AJ_0,0,0),(0,0,w')]=(0,0,AJ_0w'+(\tr AJ_0)w'), \\ 
	[(A,0,0),J(0,0,w')] &= [(A,0,0),(0,0,-J_0w')]=(0,0,-AJ_0w'-(\tr A)J_0w'),    
\end{align*} 
\noindent and
\begin{align*}
	J[J(A,0,0),J(0,0,w')] &= J[(AJ_0,0,0),(0,0,-J_0w')]\\
	&= J(0,0,Aw'-(\tr AJ_0)J_0w')\\
	&= (0,0,-J_0Aw'-(\tr AJ_0)w').
\end{align*}
\indent $(4)$ Take $X=(0,v,0)\in\a$ and $Y=(0,0,w')\in\b$. Then \eqref{eq: J es integ} holds because
\begin{align*}
	J[(0,v,0),(0,0,w')] &= J(-w'\cdot v^t,0,0)=(-w'\cdot v^tJ_0,0,0),\\
	[J(0,v,0),(0,0,w')] &= [(0,-J_0v,0),(0,0,w')]=(-w'\cdot v^t J_0,0,0),\\
	[(0,v,0),J(0,0,w')] &= [(0,v,0),(0,0,-J_0w')]=(J_0w'\cdot v^t,0,0), 
\end{align*}
\noindent and
\begin{align*}
	J[J(0,v,0),J(0,0,w')] &= J[(0,-J_0v,0),(0,0,-J_0w')]\\
	&= J(-(J_0w')\cdot(J_0v)^t,0,0)\\
	&= J(J_0w'\cdot v^tJ_0,0,0)\\
	&= (-J_0w'\cdot v^t,0,0).
\end{align*}
\indent $(5)$ Take $X,Y\in\a$ or $X,Y\in\b$. Then \eqref{eq: J es integ} holds (trivially) because these Lie subalgebras are abelian and $J$-invariant. 

\smallskip

\indent Thus, $J$ and $E$ are a complex structure and a product structure on $\sl(2n+1,\R)$, respectively. As we have observed that $J$ and $E$ anticommute, or equivalently that $J \s_{\pm} = \s_{\mp}$, we obtain that $\{J, E\}$ is a complex product structure on $\sl(2n+1,\R)$. Furthermore, the inclusion
\begin{align*}
	\gl(2n,\R) \hookrightarrow \sl(2n+1,\R)
\end{align*}
\noindent described in equation \eqref{eq: inclusion of gl in sl} is compatible with the complex product structures $\{J',E'\}$ and $\{J,E\}$. 
We summarize these facts in the following result. 
\begin{proposition} \label{prop: complex product structure}
	The Lie algebra $\sl(2n+1,\R)$ carries a complex product structure $\{J,E\}$ such that the inclusion $\gl(2n,\R) \hookrightarrow\sl(2n+1,\R)$ is compatible with the complex product structures on each Lie algebra.
\end{proposition} 
\begin{remark} 
	For $n=1$, the complex structure $J$ is equivalent to the complex structure $J_{\mathrm I_0}$ appearing in Section \ref{section: SL(3,R)}. 
\end{remark}
\indent 
Combining Propositions \ref{prop: HCX} and \ref{prop: complex product structure}, we obtain the following result. 
\begin{corollary} \label{cor: hpcx underlying}
	The real Lie algebra $\sl(2n+1,\C)_{\R}$ underlying $\sl(2n+1,\C)$ admits a hypercomplex structure such that the inclusion $\gl(2n,\C)_{\R} \hookrightarrow \sl(2n+1,\C)_{\R}$ is compatible with the hypercomplex structures on each Lie algebra. 
\end{corollary}

\indent It is immediate from Corollary \ref{cor: hpcx underlying} that the real manifold $\SL(2n+1,\C)$ carries a hypercomplex structure such that the canonical inclusion $\GL(2n,\C)\hookrightarrow \SL(2n+1,\C)$ is compatible with the left-invariant hypercomplex structures.

\begin{remark}
	The Lie group $\SL(2n+1,\C)$ admits uniform lattices, as follows from \cite[Theorem C]{Bo}. If $\Gamma$ is such a uniform lattice then $M_\Gamma:=\Gamma\backslash \SL(2n+1,\C)$ carries a hypercomplex structure naturally induced by the left-invariant one on $\SL(2n+1,\C)$. Since $\SL(2n+1,\C)$ is simply connected, we have that $\pi_1(M_\Gamma)\cong \Gamma$. Therefore, pairwise non-isomorphic uniform lattices of $\SL(2n+1,\C)$ give rise to pairwise non-homeomorphic compact hypercomplex manifolds. 
\end{remark}

\indent One may ask for a left-invariant special hyper-Hermitian metric defined on $\SL(2n+1,\C)$ that is compatible with the hypercomplex structure alluded to in Corollary \ref{cor: hpcx underlying}, a possibility of which is a HKT metric. We close this section by proving that there is no left-invariant HKT metric on $\SL(2n+1,\C)$ compatible with the left-invariant hypercomplex structure induced by the one described in Corollary \ref{cor: hpcx underlying}. As usual, it is enough to prove the non-existence at the Lie algebra level. The proof is grounded in the fact that for arbitrary $n \in \N$ we can find natural inclusions of Lie algebras
\begin{align} \label{eq: invariant inclusions}
	\sl(2n+1,\C)\hookrightarrow \sl(2n+3,\C)
\end{align}
\noindent that are compatible with the hypercomplex structures on each Lie algebra given in Corollary \ref{cor: hpcx underlying}. This inclusion is given by the map
\begin{gather*}
	\sl(2n+1,\C)\hookrightarrow \sl(2n+3,\C), \quad 
	\left[\begin{array}{cc|c} A & C & w_1 \cr B & D & w_2 \cr \hline v_1^t & v_2^t & -\tr(A+D)\cr \end{array}\right] \mapsto \left[\begin{array}{cc|c} A^\dagger & C^\dagger & w_1^\dagger \cr B^\dagger & D^\dagger & w_2^\dagger \cr \hline (v_1^\dagger)^t & (v_2^\dagger)^t & -\tr(A+D)\cr \end{array}\right],
\end{gather*}
\noindent where, if $X$ is a $n\times n$ matrix, $X^\dagger$ is the $(n+1)\times (n+1)$ matrix obtained adding one row of 0's in the bottom and a column of 0's in the right, and $w^\dagger :=\left(\begin{smallarray}{c} w \\ 0 \end{smallarray}\right)$ for any $w\in\C^n$. Note that this inclusion is compatible with the hypercomplex structure on each Lie algebra.

	\begin{proposition} \label{prop: there is no compatible HKT metric}
		There is no HKT metric defined on $\sl(2n+1,\C)$ compatible with the hypercomplex structure given in Corollary \ref{cor: hpcx underlying}.
	\end{proposition}
	\begin{proof}
		
		It follows from \eqref{eq: invariant inclusions} and Corollary \ref{cor: hpcx underlying} that we have the following sequence of inclusions compatible with the hypercomplex structures:
		\[\gl(2,\C) \hookrightarrow \sl(3,\C)\hookrightarrow \sl(5,\C)\hookrightarrow \sl(7,\C) \hookrightarrow \cdots \]
		Now, as a consequence of this chain of inclusions and Lemma \ref{lem: HKT subalgebra}, it suffices to prove that there is no HKT metric on $\gl(2,\C)$ compatible with the hypercomplex structure $\{J_1, J_2, J_3\}$ given in Corollary \ref{cor: hpcx underlying} and equation \eqref{eq: J1y2}. We will see this by means of contradiction. Let us denote
		\begin{gather*}
			e_1 :=
			\left[\begin{smallmatrix}
				1&0\\0&0
			\end{smallmatrix}\right], \quad
			e_2 := 
			\left[\begin{smallmatrix} 
				0&1\\0&0
			\end{smallmatrix}\right], \quad
			e_3 := 
			\left[\begin{smallmatrix}
				0&0\\1&0
			\end{smallmatrix}\right], \quad
			e_4 := 
			\left[\begin{smallmatrix}
				0&0\\0&1
			\end{smallmatrix}\right], \\
			e_5 := i e_1, \quad e_6 := i e_2, \quad e_7 := i e_3, \quad e_8 := i e_4. 
		\end{gather*}
		\noindent The hypercomplex structure $\{J_1, J_2\}$ is determined by the following relations: 
		\begin{gather*}
			J_1 e_1 = -e_2, \quad J_1 e_3 = -e_4, \quad J_1 e_5 = -e_6, \quad J_1 e_7 = -e_8, \\
			J_2 e_1 = \phantom{+} e_5, \quad J_2 e_2 = -e_6, \quad J_2 e_3 = \phantom{+} e_7, \quad J_2 e_4 = -e_8.
		\end{gather*} 
		\noindent Assume that there is a hyper-Hermitian metric $g$ that is HKT. Since $J_1, J_2$ and $J_3$ are orthogonal with respect to $g$, the matrix $[A]$ defined by $g$ via $[A]_{ij}=g(e_i,e_j)$ has the following form:
		\begin{align*}
			[A] := \left[\begin{array}{cccc|cccc}
				a_{11}&0&a_{13}&a_{14}&0&0&a_{17}&a_{18}\\ 
				0&a_{11}&-a_{14}&a_{13}&0&0&-a_{18}&a_{17} \\ 
				a_{13}&-a_{14}&a_{33}&0&-a_{17}&-a_{18}&0&0\\a_{14}&a_{13}&0&a_{33}&a_{18}&-a_{17}&0&0\\ \hline 0&0& -a_{17} & a_{18} &a_{11} & 0&a_{13}&-a_{14} \\ 0&0&-a_{18}&-a_{17}&0&a_{11}&a_{14}&a_{13}\\
				a_{17}&-a_{18}&0&0&a_{13}&a_{14}&a_{33}&0\\
				a_{18}&a_{17}&0&0&-a_{14}&a_{13}&0&a_{33}
			\end{array}\right].
		\end{align*}
		\indent 
		Filling in equation \eqref{eq: HKT} with the triples $(x, y, z)$ given by
		\begin{align*}
			(e_1,e_2,e_3), \quad (e_1,e_2,e_4), \quad (e_1, e_2, e_8), \quad (e_2,e_3,e_4),
		\end{align*}    
		\noindent one finds that $a_{11}=0$, $a_{13}=0$, $a_{17}=0$, $a_{33}=0$. One can then compute the eigenvalues of $[A]$ to be $\{\pm \sqrt{a_{14}^2+a_{18}^2}\}$, which contradicts the positive definiteness of the metric $g$. 
	\end{proof}
	
	\begin{corollary}
		Let $\Gamma$ be a uniform lattice of $\mathrm{SL}(2n+1,\C)$. Then the compact manifold $M_\Gamma=\Gamma\backslash \mathrm{SL}(2n+1,\C)$ does not admit any HKT metric compatible with the hypercomplex structure given in Corollary \ref{cor: hpcx underlying}.
	\end{corollary}
	
	\begin{proof}
		The statement follows from Proposition \ref{prop: there is no compatible HKT metric} and the fact that the existence of an HKT metric on $M_\Gamma$ implies the existence of an invariant one (\cite[Theorem 3.1]{FG}). 
	\end{proof}
	
	\begin{remark}
		Proposition \ref{prop: there is no compatible HKT metric} also holds for $\gl(2n,\C)$ equipped with the hypercomplex structure given in \cite{AS} also obtained  using Proposition \ref{prop: HCX}. Indeed, there is another chain of hypercomplex inclusions given by
		\[ \gl(2,\C)\hookrightarrow \gl(4,\C)\hookrightarrow \gl(6,\C) \hookrightarrow \cdots\!,\] where each inclusion is given by the map \begin{gather*}
			\gl(2n,\C)\hookrightarrow \gl(2n+2,\C), \quad \begin{bmatrix}
				A & C \\ B & D 
			\end{bmatrix} \mapsto \begin{bmatrix}
				A^\dagger & C^\dagger \\ B^\dagger & D^\dagger 
			\end{bmatrix}.
		\end{gather*}
		Since we already know that $\gl(2,\C)$ does not admit an HKT metric compatible with the induced hypercomplex structure, it follows again from Lemma \ref{lem: HKT subalgebra} that $\gl(2n,\C)$ does not admit an HKT metric compatible with the given hypercomplex structure.  
	\end{remark} 
	
	\medskip 
	
	\section{Holonomy of the Obata connection on \texorpdfstring{$\SL(2n+1,\C)$}{}} \label{sec: holonomy}
	
	\indent This last section is devoted to computing the holonomy of the Obata connection $\nob$ of the hypercomplex structure on $\SL(2n+1,\C)$ obtained in Corollary \ref{cor: hpcx underlying} from a complex product structure $\{J, E\}$ on $\SL(2n+1,\R)$. Recall from Proposition \ref{prop: complex-hol} that the holonomy of $\nob$ is the complexification of the holonomy of  the connection $\ncp$ associated to the complex structure $\{J,E\}$. As a first step, we apply Proposition \ref{prop: existe ncp} or, equivalently, equation \eqref{eqs: ncp-casos} to compute $\ncp$ explicitly. 
	\begin{proposition}\label{prop: nablacp}
		Let $x,y\in \sl(2n+1,\R)$, given by \[x=\left[\begin{array}{cc|c} A&C&w_1\\B&D&w_2\\ \hline v_1^t&v_2^t&-\tr (A+D)\end{array}\right], \quad y=\left[\begin{array}{cc|c}E&G&z_1\\F&H&z_2\\ \hline y_1^t&y_2^t&-\tr (E+H)\end{array}\right],\] and let $x_{\pm}:=\pi_{\pm} x$ and $y_{\pm}:=\pi_{\pm}y$, with $\pi_{\pm}$ the projections to $\mathfrak{s}_+$ and $\mathfrak{s}_-$, respectively. Then,
		\begin{enumerate} 
			\item[$\ri$] $\ncp_{x_+} y_+=\left[\begin{array}{cc|c} AE&\boldsymbol{0}&\boldsymbol{0}\\BE+w_2 y_1^t&\boldsymbol{0}&Az_2+Ew_2+\tr (A)z_2 \\ \hline v_1^t E-\tr(A)y_1^t &\boldsymbol{0}&-\tr(AE)\end{array}\right]$.
			\smallskip
			\item[$\rii$] $\ncp_{x_-} y_+=\left[\begin{array}{cc|c}CF+w_1 y_1^t&\boldsymbol{0}&\boldsymbol{0}\\DF&\boldsymbol{0}&Dz_2-Fw_1+\tr(D) z_2\\ \hline v_2^t F-\tr(D) y_1^t  &\boldsymbol{0} &-\tr(CF+w_1 y_1^t)\end{array}\right]$.
		\end{enumerate}
		The values of $\ncp_{x_+} y_-$ and $\ncp_{x_-} y_-$ can be computed from the expressions above and $\ncp J=0$.
	\end{proposition}
	
	\smallskip
	
	\indent Proposition \ref{prop: nablacp} comes in handy to compute the curvature operators and their covariant derivatives. Proposition \ref{prop: existe ncp} ensures that $R^{\operatorname{CP}}(\s_+,\s_+)\equiv 0$ and $R^{\operatorname{CP}}(\s_-,\s_-)\equiv 0$, thus we only need to consider the curvature endomorphisms $R^{\operatorname{CP}}(x_+,y_-)$ for $x_+\in\s_+$, $y_-\in\s_-$. The following result can be obtained by tedious yet straightforward calculations.  
	
	\begin{proposition}\label{prop: curv nablacp}
		Let $x,y$ be as in Proposition \ref{prop: nablacp}, and let $z_+=\left[\begin{array}{cc|c} U&\boldsymbol{0}&\boldsymbol{0}\\V&\boldsymbol{0}&q_2\\ \hline p_1^t&\boldsymbol{0}&-\tr U\end{array}\right]\in \s_+$. 
		
		\noindent Let $\xi(x_+,y_-,z_+)\in \R^n$ be the following expression: 
		\begin{align*}
			\xi(x_+,y_-,z_+) &:= ([A,H]+[G,B])q_2+([V,A]+[B,U])z_1+([G,V]+[U,H])w_2\\
			&\quad+(\langle v_1,z_1\rangle-\langle y_2,w_2\rangle )q_2 + (\langle v_1,q_2\rangle +\langle w_2,p_1\rangle) z_1 +(\langle p_1,z_1\rangle -\langle y_2,q_2\rangle) w_2.
		\end{align*}
		Then $R^{\operatorname{CP}}(x_+,y_-)z_+= \left[\begin{array}{cc|c}\boldsymbol{0}&\boldsymbol{0}&\boldsymbol{0}\\\boldsymbol{0}&\boldsymbol{0}& \xi(x_+,y_-,z_+)\\ \hline \boldsymbol{0}&\boldsymbol{0}&0\end{array}\right].$
		The values of $R^{\operatorname{CP}}(x_+,y_-)z_-$ for $z_-\in \s_-$ can be computed from this expression and $R^{\mathrm{CP}}(x_+,y_-) J= J R^{\mathrm{CP}}(x_+,y_-)$.
	\end{proposition}
	
	We state some particular cases which will be used in the proof of Theorem \ref{thm: holonomy}:
	\begin{gather}
		x_+=v_1^t,\, y_-=z_1, \, z_+=q_2: \quad \xi(x_+,y_-,z_+) = \la v_1, q_2\ra z_1+\la v_1, z_1\ra q_2,\label{eq: curv interest}\\
		x_+=w_2, \, y_-=z_1, \, z_+=p_1:\quad \xi(x_+,y_-,z_+)=\la p_1, w_2\ra z_1+\la p_1, z_1\ra w_2. \label{eq: curv interest2}
	\end{gather}
	In the expressions above, when we write for instance $x_+=w_2$, we mean that the matrix $x_+\in\s_+$ as in Proposition \ref{prop: nablacp} has all its components equal to zero, except possibly by the vector $w_2$. 
	
	\medskip
	
	\noindent We compute next the covariant derivatives of the curvature endomorphisms.
	
	\begin{proposition}\label{prop: dercov nablacp}
		Let $x,y,z_+\in \sl(2n+1,\R)$ be as in Proposition \ref{prop: curv nablacp}. Let also \[u_+=\left[\begin{array}{cc|c} X&\boldsymbol{0}&\boldsymbol{0}\\Y&\boldsymbol{0}&s_2\\ \hline r_1^t&\boldsymbol{0}&-\tr X\end{array}\right] \quad \text{and}\quad u_-=\left[\begin{array}{cc|c}\boldsymbol{0}&Z&s_1\\\boldsymbol{0}&W&\boldsymbol{0}\\ \hline \boldsymbol{0}&r_2^t&-\tr W\end{array}\right].\]
		Then 
		\begin{enumerate} \item[$\ri$] $[\ncp_{u_+} R^{\operatorname{CP}}(x_+,y_-)]z_+=\left[\begin{array}{cc|c}
				\boldsymbol{0}&\boldsymbol{0}&\boldsymbol{0}\\\boldsymbol{0}&\boldsymbol{0}&\nu_+(u_+,x_+,y_-,z_+)\\ \hline \boldsymbol{0}&\boldsymbol{0}&0
			\end{array}\right]$, where $\nu_+(u_+,x_+,y_-,z_+)\in \R^n$ is 
			\begin{align*}
				\nu_+(u_+,x_+,y_-,z_+)&=([A,YU]+X[V,A]+[X,B]U+(\la w_2, p_1\ra+\la q_2, v_1\ra )X)z_1\\
				&+(-\la Us_2, v_1\ra-\la X q_2, v_1\ra-\la Uw_2, r_1\ra)z_1\\
				&+([YU,G]+X[G,V]+[H,X]U+(\la z_1, p_1\ra -\la q_2, y_2\ra) X)w_2\\
				&+(\la Us_2, y_2\ra+\la Xq_2, y_2\ra-\la Uz_1, r_1\ra)w_2\\
				&+[X,[A,H]+[G,B]]q_2\\
				&+([H,A]U+[B,G]U+\la y_2, w_2\ra U-\la v_1,z_1\ra U) s_2\\
				&+(\la G w_2, p_1\ra-\la Az_1, p_1\ra+\la z_1, p_1\ra A-\la w_2, p_1\ra G) s_2\\
				&+\tr(X)( ([G,V]+[U,H]+2\la z_1, p_1\ra)w_2+([V,A]+[B,U]+2\la w_2, p_1\ra)z_1).
			\end{align*}
			\item[$\rii$] $[\ncp_{u_-} R^{\operatorname{CP}}(x_+,y_-)]z_+=\left[\begin{array}{cc|c}
				\boldsymbol{0}&\boldsymbol{0}&\boldsymbol{0}\\\boldsymbol{0}&\boldsymbol{0}&\nu_-(u_-,x_+,y_-,z_+)\\ \hline \boldsymbol{0}&\boldsymbol{0}&0
			\end{array}\right]$,
			where $\nu_-(u_+,x_+,y_-,z_+)\in \R^n$ is 
			\begin{align*}
				\nu_-(u_-,x_+,y_-,z_+)&=([ZV,B]+[A,W]V+W[B,U]+(\la v_1,q_2\ra+\la p_1,w_2\ra)W)z_1\\
				&\quad+(\la v_1, Vs_1\ra-\la v_1, Wq_2\ra-\la r_2, Vw_2\ra) z_1\\
				&\quad+ ([H,ZV]+[W,G]V+W[U,H]+(\la p_1,z_1\ra-\la y_2, q_2\ra)W)w_2\\
				&\quad+(\la y_2,Wq_2\ra -\la y_2, Vs_1\ra-\la r_2, Vs_1\ra)w_2\\
				&\quad+[W,[A,H]+[G,B]]q_2\\
				&\quad+(([A,H]+[G,B])V-\la y_2, w_2\ra V+\la v_1, z_1\ra V-\la p_1, z_1\ra B+\la p_1,w_2\ra H)s_1\\
				&\quad+(\la p_1, Bz_1\ra-\la p_1,Hw_2\ra)s_1\\
				&\quad+\tr(W)\big(([U,H]+[G,V]+2\la p_1, z_1\ra \Id)w_2+([B,U]+[V,A]+2\la p_1, w_2\ra\Id)z_1\big).
			\end{align*}
		\end{enumerate}
		Moreover,  $[\ncp_{u_\pm } R^{\operatorname{CP}}(x_+,y_-)]z_-$ for $z_-\in\s_-$ can be computed from  the expressions above and the fact that $[\ncp_{u_\pm} R^{\operatorname{CP}}(x_+,y_-)]J=J[\ncp_{u_\pm} R^{\operatorname{CP}}(x_+,y_-)]$.
	\end{proposition}
	
	In particular, we obtain the following expressions that will be used in the proof of Theorem \ref{thm: holonomy}:
	\begin{gather}
		u_+=r_1^t, \, x_+=w_2, y_-=z_1,\, z_+=U:\quad \nu_+(u_+,x_+,y_-,z_+)=-\la r_1, Uw_2\ra z_1-\la r_1,Uz_1\ra w_2, \label{eq: dercov interest 1}\\ 
		u_-=r_2^t, \, x_+=w_2, \, y_-=z_1, z_+=V:\quad \nu_-(u_-,x_+,y_-,z_+)=-\la r_2, Vz_1\ra w_2-\la r_2, Vw_2\ra z_1. \label{eq: dercov interest 2}
	\end{gather}

	\medskip
	
	Using Theorem \ref{thm: hol} and the expressions for the curvature operators and their covariant derivatives obtained in Propositions \ref{prop: curv nablacp} and \ref{prop: dercov nablacp} we determine the holonomy algebra $\hol^{\mathrm{CP}}$ of the connection $\ncp$ on $\sl(2n+1,\R)$. We refer to the decomposition given in equation \eqref{eq: decomposition of sl in gl a b}. 
	
	\begin{theorem}\label{thm: holonomy}
		The holonomy algebra $\hol^{\mathrm{CP}}$ of the connection $\ncp$ on $\sl(2n+1,\R)$ is given by
		\begin{equation}\label{eq: hol-CP}
			\hol^{\mathrm{CP}}=\{ T\in \operatorname{End}(\sl(2n+1,\R)) : [T,J]=0, \; [T,E]=0, \; \operatorname{Im}T\subseteq \b \}.
		\end{equation}  
	\end{theorem}
	
	\begin{proof} 
		In this proof we will use the notation $\s:=\sl(2n+1,\R)$, for simplicity. 
		Let us denote $\h$ the right-hand side in \eqref{eq: hol-CP}, which is a Lie subalgebra of $\operatorname{End}(\s)$.
		
		Due to $\ncp J=0$, $\ncp E=0$ and Proposition \ref{prop: curv nablacp}, 
		it is clear that $\hol^{\mathrm{CP}}\subseteq \h$. We will show next the reverse inclusion.
		
		First we note that if $T\in \h$ then both $\s_+$ and $\s_-$ are $T$-invariant. Let us set accordingly $T_+= T|_{\s_+}:\s_+\to \s_+$ and $T_-= T|_{\s_-}:\s_-\to \s_-$. Moreover, if $\mathcal{B}_+$ is a basis of $\s_+$ then the set $\mathcal{B}_-$ obtained by applying $J$ to the elements of $\mathcal{B}_+$ is a basis of $\s_-$. Since $T$ commutes with $J$ we have that $[T_+]_{\mathcal{B}_+}=[T_-]_{\mathcal{B}_-}$, hence the matrix of $T$ in the basis $\mathcal{B}=\mathcal{B}_+\cup \mathcal{B}_-$ of $\s$ is given by $[T]_{\mathcal B}=A^{\oplus 2}$, where $A=[T_+]_{\mathcal{B}_+}=[T_-]_{\mathcal{B}_-}$.
		
		Next, we will show that a certain subspace $\h_1$ of $\h$ is contained in $\mathfrak{hol}^{\mathrm{CP}}$. Indeed, we will determine a family of endomorphisms $T_+:\s_+\to \s_+$ with $\operatorname{Im}T_+ \subseteq \b_+$ that can be obtained as a linear combination of curvature endomorphisms and their covariant derivatives, restricted to $\s_+$. In order to do so, we describe more precisely the basis $\mathcal{B}_+$ of $\s_+$. Let us set
		\[
		\alpha_{ij}=\left[\begin{array}{cc|c} E_{ij} & \boldsymbol{0} &\boldsymbol{0} \cr \boldsymbol{0} & \boldsymbol{0} & \boldsymbol{0} \cr  \hline\boldsymbol{0} & \boldsymbol{0} & -\delta_{ij}\cr \end{array}\right], \quad \beta_{ij}=\left[\begin{array}{cc|c} \boldsymbol{0} & \boldsymbol{0} &\boldsymbol{0} \cr E_{ij} & \boldsymbol{0} & \boldsymbol{0} \cr \hline \boldsymbol{0} & \boldsymbol{0} & 0\cr \end{array}\right], \quad \gamma_{i}=\left[\begin{array}{cc|c} \boldsymbol{0} & \boldsymbol{0} &\boldsymbol{0} \cr \boldsymbol{0} & \boldsymbol{0} & \boldsymbol{0} \cr \hline e_i^t & \boldsymbol{0} & 0\cr \end{array}\right], \quad \mu_{i}=\left[\begin{array}{cc|c} \boldsymbol{0} & \boldsymbol{0} &\boldsymbol{0} \cr \boldsymbol{0} & \boldsymbol{0} & e_i \cr \hline\boldsymbol{0} & \boldsymbol{0} & 0\cr \end{array}\right],
		\]
		where $E_{ij}$ is the usual $n\times n$ matrix whose only non-zero entry is equal to $1$ in the spot $(i,j)$, and $\{e_i\}$ is the canonical basis of $\R^n$. Therefore, 
		\[ \mathcal{B}_+=\{ \alpha_{ij}, \, \beta_{ij} : 1\leq i,j\leq n\} \cup \{\gamma_i, \, \mu_i : 1\leq i\leq n\}.\]
		Note that $\{ \alpha_{ij}, \, \beta_{ij} : 1\leq i,j\leq n\}$ is a basis of $\gl(n,\R)\cap \s_+$, $\{\gamma_i : 1\leq i\leq n\}$ is a basis of $\a_+$ and $\{\mu_i : 1\leq i\leq n\}$ is a basis of $\b_+$. Let us also denote $\l_1:=\text{span}\{\alpha_{ij}\}$ and $\l_2:=\text{span}\{\beta_{ij}\}$, so that $\gl(n,\R)\cap \s_+=\l_1\oplus \l_2$.
		
		We divide our analysis in several stages.
		
		\noindent (1) We show here that given $\mu_i$ and $\mu_j$, with $1\leq i,j\leq n$, there is an operator in $\hol^{\mathrm{CP}}$ that sends $\mu_i$ to $\mu_j$, and $\mu_k$ to $0$, where $k\neq i$. Indeed, let us consider 
		\[ x_+=\gamma_i, \quad y_-=J\mu_j=\left[\begin{array}{cc|c} \boldsymbol{0} & \boldsymbol{0} &e_j \cr \boldsymbol{0} & \boldsymbol{0} & \boldsymbol{0} \cr \hline \boldsymbol{0} & \boldsymbol{0} & 0\cr \end{array}\right].\]
		Then, according to \eqref{eq: curv interest}, we have that 
		\begin{equation}\label{eq: prueba holonomy} R^{\operatorname{CP}}(\gamma_i,J\mu_j)\mu_k = \langle e_i, e_k\rangle \mu_j+ \langle e_i, e_j\rangle \mu_k.
		\end{equation}
		Hence, if $i\neq j$ we have $R^{\operatorname{CP}}(\gamma_i,J\mu_j)\mu_k=\delta_{ik} \mu_j$ and we are done in this case. On the other hand, if $i=j$, let us consider $S_i\in \hol^{\mathrm{CP}}$ given by
		\[ S_i=\frac{1}{n-1}\left[(n-2)R^{\operatorname{CP}}(\gamma_i,J\mu_i)-\sum_{k\neq i}R^{\operatorname{CP}}(\gamma_k,J\mu_k)\right].\]
		Then, using that $R^{\operatorname{CP}}(\gamma_r,J\mu_r)\mu_s=(1+\delta_{rs})\mu_s$ (which follows from \eqref{eq: prueba holonomy}), it is easily verified that $S_i\mu_k=\delta_{ik}\mu_k$.
		
		Moreover, $R^{\operatorname{CP}}(\gamma_i,J\mu_j)$ vanishes on $\l_1\oplus\l_2\oplus\a_+$, for any $i,j$, due to Proposition \ref{prop: curv nablacp}. As a consequence, all the operators $T_+:\b_+\to \b_+$ (that is, $T_+:\s_+\to \b_+$ satisfying $T|_{\l_1\oplus\l_2\oplus\a_+}=0$) are a linear combination of curvature endomorphisms restricted to $\s_+$, as we wanted to show.
		
		\smallskip
		
		\noindent (2) We show next that any operator $T_+:\a_+\to \b_+$ (that is, $T_+:\s_+\to \b_+$ with $T|_{\l_1\oplus\l_2\oplus\b_+}=0$) whose matrix in the bases $\{\gamma_i\}$ of $\a_+$ and $\{\mu_i\}$ of $\b_+$ is \textit{symmetric} can be obtained as a linear combination of curvature endomorphisms restricted to $\s_+$. Indeed, for any $1\leq i,j\leq n$, we have that $R^{\operatorname{CP}}(\mu_i,J\mu_j)|_{\l_1\oplus\l_2\oplus\b_+}=0$ (due to Proposition \ref{prop: curv nablacp}) and on $\a_+$ we compute (using \eqref{eq: curv interest2})
		\[ R^{\operatorname{CP}}(\mu_i,J\mu_j)\gamma_k= \langle e_i,e_k\rangle \mu_j+\langle e_j, e_k\rangle \mu_i =\delta_{ik}\mu_j+\delta_{jk}\mu_i.  \]
		Therefore, the matrix of $R^{\operatorname{CP}}(\mu_i,J\mu_j):\a_+\to\b_+$ is $E_{ij}+E_{ji}$, and thus any $T_+$ as above can be written as a linear combination of these operators. 
		
		\smallskip
		
		\noindent (3) We prove now that certain operators $T_+:\l_1\to \b_+$ (that is, $T_+:\s_+\to \b_+$ with $T|_{\l_2\oplus\a_+\oplus\b_+}=0$) are a linear combination of covariant derivatives of the curvature endomorphisms, hence they are in $\hol^{\mathrm{CP}}$. For any $1\leq i,j,k\leq n$, we have that $\ncp_{\gamma_i}R^{\operatorname{CP}}(\mu_j, J\mu_k)|_{\l_2\oplus\a_+\oplus\b_+}=0$ (as follows from Proposition \ref{prop: dercov nablacp}) and  on $\l_1$ we compute, using equation \eqref{eq: dercov interest 1},
		\[ [\ncp_{\gamma_i}R^{\operatorname{CP}}(\mu_j, J\mu_k)] \alpha_{rs}= -\langle e_i, E_{rs} e_j\rangle \mu_k - \langle e_i, E_{rs} e_k\rangle \mu_j=-\delta_{js}\delta_{ir}\mu_k-\delta_{ks}\delta_{ir}\mu_j.\]
		Hence, if $r\neq i$ this expression vanishes. However, if $r=i$ we have 
		\[ [\ncp_{\gamma_i}R^{\operatorname{CP}}(\mu_j, J\mu_k)] \alpha_{is}= -\delta_{js}\mu_k- \delta_{ks}\mu_j.  \]
		That is, $\ncp_{\gamma_i}R^{\operatorname{CP}}(\mu_j, J\mu_k)$ sends $\alpha_{ij}$ to $-\alpha_k$ and $\alpha_{ik}$ to $-\alpha_j$, and vanishes on any other $\alpha_{rs}$. Setting $\l_1^i=\text{span}\{\alpha_{i1},\ldots,\alpha_{in} \}$, this means that the matrix of $\ncp_{\gamma_i}R^{\operatorname{CP}}(\mu_j, J\mu_k)|_{\l_1^i}:\l_1^i\to \b_+$ is given by $-E_{jk}-E_{kj}$. Varying $j$ and $k$, we obtain that if $T_+:\l_1\to \b_+$ is an operator whose only non-vanishing component is $T_+|_{\l_1^i}$ and its matrix in the bases $\{\alpha_{i1},\ldots,\alpha_{in}\}$ and $\{\mu_1,\ldots,\mu_n\}$ is \textit{symmetric} then $T_+$ is in $\mathfrak{hol}^{\mathrm{CP}}$.
		
		\smallskip
		
		\noindent (4) Consider now an operator $T_+:\l_2\to \b_+$ (that is, $T_+:\s_+\to \b_+$ with $T|_{\l_1\oplus\a_+\oplus\b_+}=0$), and set $\l_2^i=\text{span}\{\beta_{i1},\ldots,\beta_{in} \}$. In a similar manner to the one in (3), we can see that if the only non-vanishing component of $T_+$ is $T_+|_{\l_2^i}$ and its matrix in the bases $\{\beta_{i1},\ldots,\beta_{in}\}$ and $\{\mu_1,\ldots,\mu_n\}$ is \textit{symmetric} then $T_+$ is in $\hol^{\mathrm{CP}}$. Here we use equation \eqref{eq: dercov interest 2}.
		
		\medskip
		
		To sum up, we have shown that the following subspace $\h_1$ of $\h$ is contained in $\hol^{\mathrm{CP}}$. Here we are identifying an operator $T\in \h$ with its matrix using the ordered basis $\{\alpha_{ij}\}\cup\{\beta_{ij}\}\cup \{\gamma_i\}\cup \{\mu_i\}$, using the lexicographical order in the subscripts:
		\[ \h_1=\left\{M_{(A_1,\ldots,A_n,B_1,\ldots,B_n,C,D)}\mid A_i,B_i,C,D\in \gl(n,\R), \, A_i^t=A_i, B_i^t=B_i, C^t=C \, \forall i \right\}, \]
		where 
		\begin{equation}\label{eq: matrixM}
			M_{(A_1,\ldots,A_n,B_1,\ldots,B_n,C,D)}=
			\left[\begin{array}{ccc|ccc|c|c}
				0 & \cdots & 0  & 0 &\cdots & 0 & 0 & 0 \\
				\vdots  & & \vdots  & \vdots  & & \vdots  & \vdots&\vdots \\
				0 & \cdots & 0  & 0 &\cdots & 0 & 0 & 0 \\
				\hline 
				A_1 & \cdots & A_n & B_1 & \cdots & B_n & C & D
			\end{array}\right]^{\oplus 2}. 
		\end{equation}
		Since $\hol^{\mathrm{CP}}$ is a Lie subalgebra of endomorphisms, we have that the Lie bracket of two elements in $\h_1$ is in $\hol^{\mathrm{CP}}$. A generic element in $\h$ is a matrix as in \eqref{eq: matrixM}, but without asking $A_i,B_i,C$ to be symmetric. Note that 
		\begin{align}\label{eq: bracketM}
			[M_{(A_1,\ldots,A_n,B_1,\ldots,B_n,C,D)}, \, & M_{(A_1',\ldots,A_n',B_1',\ldots,B_n',C',D')}]= \\ & =M_{(DA_1'-D'A_1,\ldots,DA_n'-D'A_n,DB_1'-D'B_1,\ldots,DB_n'-D'B_n,DC'-D'C,[D,D'])}. \nonumber
		\end{align}
		If $X\in\gl(n,\R)$ is an arbitrary matrix, let us set in \eqref{eq: bracketM} the following: in the first matrix, everything is $0$ except for $D=X$, and in the second matrix everything is $0$ except for $A_i'=I_n$, the $n\times n$ identity matrix. Then the resulting commutator has every entry equal to $0$, except for the spot $A_i=X$, and it is in the holonomy algebra. Repeating this for the spots $B_i$ and $C$, we obtain that a matrix $M_{(A_1,\ldots,A_n,B_1,\ldots,B_n,C,D)}$ with arbitrary $n\times n$ matrices $A_i, B_i, C$ and $D$ is in $\hol^{\mathrm{CP}}$. This shows that $\h\subseteq \hol^{\mathrm{CP}}$, and the proof is complete.
	\end{proof}
	
	\medskip
	
	Let us denote by $\mu:\gl(n,\R)\to \operatorname{End}(\gl(n,\R))$ the representation of $\gl(n,\R)$ on $\gl(n,\R)$ given by left multiplication of matrices: $\mu(A)X=AX$. For any $k\in \N$, let $(V_k,\mu_k)$ denote the direct sum representation $(V_k,\mu_k)=(\gl(n,\R),\mu)\oplus \cdots \oplus (\gl(n,\R),\mu)$ ($k$ times). The following result is a consequence of the proof of Theorem \ref{thm: holonomy} (see \eqref{eq: bracketM}).
	
	\begin{corollary}\label{coro: hol-CP}
		The holonomy algebra $\mathfrak{hol}^{\mathrm{CP}}$ corresponding to the complex product structure $\{J,E\}$ on $\sl(2n+1,\R)$ is isomorphic to 
		\[ \hol^{\mathrm{CP}}\cong \gl(n,\R)\ltimes_{\mu_{2n+1}} V_{2n+1}. \] 
		The abelian ideal $V_{2n+1}$ and the Lie algebra $\hol^{\mathrm{CP}}$ have dimension $n^2(2n+1)$ and $n^2(2n+2)$, respectively.
	\end{corollary}

	
	We are finally in position to determine the holonomy of the Obata connection on $\sl(2n+1,\C)_\R$.
	
	\begin{theorem}
		The holonomy algebra $\mathfrak{hol}^{\mathrm{Ob}}$ corresponding to the hypercomplex structure $\hcx_{\alpha=1,2,3}$ on $\sl(2n+1,\C)$ is isomorphic to the Lie algebra 
		\[ \mathfrak{hol}^{\mathrm{Ob}}\cong \gl(n,\C) \ltimes_{\mu^\C_{2n+1}} V^\C_{2n+1}.\] 
		The abelian ideal $V^\C_{2n+1}$ and $\hol^{\mathrm{Ob}}$ have real dimension $2n^2(2n+1)$  and  $2n^2(2n+2)$, respectively.
	\end{theorem}
	
	\begin{proof}
		This result follows in a straightforward manner from Corollary \ref{coro: hol-CP} and Proposition \ref{prop: complex-hol}.
	\end{proof}
	
	\begin{remark}
		For $n=1$ this Obata holonomy algebra is solvable, and for $n>1$, the semisimple factor in the Levi decomposition is $\sl(n,\C)$.
	\end{remark}
	
	
	\smallskip
	
	Let us set $d_n:=\dim_\R \sl(2n+1,\C)=2\dim \sl(2n+1,\R)=8n(n+1)$ and, for $\mathbb{K}=\R$ or $\C$ and $m,t\in \N$ with $m>t$,   \[\mathfrak{v}(m,t,\mathbb{K}):=\{A\in \gl(m,\mathbb{K}): A_{ij}=0\; \text{for all } i\leq m-t\}.\]
	As matrix Lie algebras, we may regard the holonomy Lie algebras we have been working with as follows:
	\begin{align*}
		\hol^{\mathrm{CP}}&=\left\{ A^{\oplus 2}\in \gl(\tfrac12 d_n,\R):  A\in \mathfrak{v}(\tfrac14 d_n, n, \R)\right\}, \\  
		\hol^{\mathrm{Ob}}&=\left\{ A^{\oplus 2}\in \gl(\tfrac12 d_n,\C):  A\in \mathfrak{v}(\tfrac14 d_n ,n,\mathbb{C})\right\}.
	\end{align*}
	In terms of the real representation of quaternionic matrices given in \eqref{eq: H} it is clear that
	\[ \hol^{\mathrm{Ob}}=\left\{ H(X,Y,0,0)\in \gl(\tfrac14 d_n,\H)\subseteq \gl(d_n,\R) : X,Y\in \mathfrak{v}(\tfrac14 d_n, n,\R) \right\}.\]

	\begin{corollary}\label{coro: hol en SL}
		The Obata holonomy of $\SL(2n+1,\C)$ is not contained in $\SL(\tfrac14 d_n,\mathbb{H})$.
	\end{corollary}
	\begin{proof}
		This follows from Corollary \ref{coro: holonomies} and the fact that $\mathfrak{v}(\tfrac14 d_n, n,\R)$ is not contained in $\sl(\tfrac14 d_n,\R)$.
	\end{proof}
	
	Concerning the Obata holonomy group $\operatorname{Hol}^{\mathrm{Ob}}$ associated to the hypercomplex structure given on $\SL(2n+1,\C)$, let us note that it coincides with the restricted holonomy group $\operatorname{Hol}_0^{\mathrm{Ob}}$ since $\SL(2n+1,\C)$ is simply connected. In particular, $\operatorname{Hol}^{\mathrm{Ob}}$ is connected. It follows that $\operatorname{Hol}^{\mathrm{Ob}}$, as a complex Lie group, is the following group of complex matrices
	\[ \operatorname{Hol}^{\mathrm{Ob}}=\left\{ \left[\begin{array}{c|c} I_{\frac12 d_n-n} & \boldsymbol{0} \cr \hline X & D \cr \end{array}\right]^{\oplus 2} : D\in\GL(n,\C), \, X\in \operatorname{Mat}(n\times (2n+1)n,\C)\right\}.\]
	
	\medskip 
	
	\begin{remark}
		The canonical bundle of the compact complex manifold $(\Gamma\backslash\SL(2n+1,\C),J_\alpha)$, where $\Gamma$ is a uniform lattice of $\SL(2n+1,\C)$, is not trivial for any $\alpha=1,2,3$. Indeed, according to \cite[Theorem 5.3]{AT}, if the canonical bundle of $(\Gamma\backslash \SL(2n+1,\C),J_\alpha)$ were trivial then we would have that $\psi_\alpha\equiv 0$, where $\psi_\alpha$ is the Koszul $1$-form on $\sl(2n+1,\C)$ associated to $J_\alpha$, defined by $\psi_\alpha(x)=\tr(J_\alpha \circ \ad x)$. Moreover, by \cite[Theorem 6.1]{AT}, it is enough to see that there exists an element $\hat{x}\in \sl(2n+1,\C)$ such that $\psi_1(\hat{x})\neq 0$. Let
		\[ \hat{x}=x^{\oplus 2}, \quad \text{where} \quad x=\left[\begin{array}{cc|c} \boldsymbol{0} & \boldsymbol{0} & \boldsymbol{0}\\ E_{11} & \boldsymbol{0} & \boldsymbol{0}\\   \hline \boldsymbol{0} & \boldsymbol{0} & 0 \end{array}\right]\in \sl(2n+1,\R).\] 
		By decomposing $\sl(2n+1,\C)=(\s_+\oplus \s_-)\oplus i(\s_+\oplus\s_-)$ (where $\s_+$ and $\s_-$ are the Lie subalgebras of $\sl(2n+1,\R)$ from Section \ref{section: SL(2n+1,C)}), we have that $\psi_1 (\hat{x})=2\tr(A \circ \ad x)$, where $A=\matriz{\boldsymbol{0} & -\I_{2n(n+1)}\\ \I_{2n(n+1)}& \boldsymbol{0}}\in \End(\sl(2n+1,\R))$. However, one can compute $\tr(A\circ \ad x)=2n+2\neq 0$, which says that the canonical bundle of $(\Gamma\backslash \SL(2n+1,\C),J_\alpha)$ is not trivial for any $\alpha=1,2,3$. 
		
		This fact provides another way to show that the Obata holonomy of $(\Gamma\backslash\SL(2n+1,\C),J_\alpha)$ is not contained in the quaternionic special linear group, due to \cite[Claim 1.2]{Ver}.
	\end{remark}
	
	\begin{remark}
		There is another way to construct a hypercomplex structure on the Lie algebra $\sl(2n+1,\C)$, noting that $\sl(2n+1,\C)$ is isomorphic to $\su(2n+1)\oplus i \,\su(2n+1)$ as real Lie algebras. Since $\su(2n+1)$ admits a hypercomplex structure (\cite[Section 4.3]{Joyce}, \cite[Corollary 4.16]{DT1}), then $\sl(2n+1,\C)$ admits an induced hypercomplex structure. Its Obata holonomy is more difficult to determine than for the hypercomplex structure given in Corollary \ref{cor: hpcx underlying} and we will not pursue this here. 
	\end{remark}

	\ 
	
\end{document}